\documentstyle[epic,eepic,psfig,12pt]{article}
\font\goth=eufm10 scaled 1200

\newtheorem{example}{Example}[section]

\newtheorem{theorem}[example]{Theorem}

\def\boxit#1#2{\setbox1=\hbox{\kern#1{#2}\kern#1}%
\dimen1=\ht1 \advance\dimen1 by #1 \dimen2=\dp1 \advance\dimen2 by #1
\setbox1=\hbox{\vrule height\dimen1 depth\dimen2\box1\vrule}%
\setbox1=\vbox{\hrule\box1\hrule}%
\advance\dimen1 by .4pt \ht1=\dimen1
\advance\dimen2 by .4pt \dp1=\dimen2 \box1\relax}

\def\Proof{\medskip\noindent {\it Proof --- \ }}
\def\cqfd{\hfill $\Box$ \bigskip}
\def\adots{\mathinner{\mkern2mu\raise1pt\hbox{.}
\mkern3mu\raise4pt\hbox{.}\mkern1mu\raise7pt\hbox{.}}}
\def\<{\langle}
\def\>{\rangle}

\def\cf{{\it cf.$\ $}}
\def\ie{{\it i.e. }}

\def\S{{\bf S}}

\newfont{\bb}{cmbx10}

\def\Z{{\bf Z}}

\def\Q{{\bf Q}}

\def\Q{{\bf Q}}

\def\mod{{\rm\ mod\ }}
\def\C{{\bf C}}
\def\K{{\bf K}}
\def\O{{\cal O}}

\def\g{\hbox{\goth g}\hskip 1pt}
\def\gl{\hbox{\goth gl}\hskip 2.5pt}
\def\sl{\hbox{\goth sl}\hskip 2.5pt}
\def\sp{\hbox{\goth sp}\hskip 2.5pt}
\def\tab{{\rm Tab\, }}

\def\ker{{\rm Ker\,}}
\def\End{{\rm End\,}}
\def\e{{\bf e}}
\def\d{{\bf d}}

\title{ CRYSTAL GRAPHS AND $q$-ANALOGUES
\\ OF WEIGHT MULTIPLICITIES
\\ FOR THE ROOT SYSTEM $A_n$\thanks{Partially supported
by PRC Math-Info and EEC grant n$^0$ ERBCHRXCT930400}}

\author{Alain {\sc Lascoux}\thanks{L.I.T.P., Universit\'e
Paris 7, 2 place Jussieu, 75251 Paris cedex 05, France} ,
\rm Bernard {\sc   Leclerc}$^\dagger$ \rm and Jean-Yves
{\sc Thibon}\thanks{Institut Gaspard Monge, Universit\'e de
Marne-la-Vall\'ee, 2 rue de la Butte-Verte, 93166 Noisy-le-Grand cedex,
France}}
\date{}

\begin{document}

\maketitle

\begin{abstract}
We give an expression of the $q$-analogues of the multiplicities of
weights in irreducible $\sl_{n+1}$-modules in terms of the geometry
of the crystal graph attached to the corresponding $U_q(\sl_{n+1})$-modules.
As an application, we describe multivariate polynomial analogues
of the multiplicities of the zero weight, refining Kostant's generalized
exponents.
\end{abstract}


\section{Introduction}

There exist interesting $q$-analogues of the multiplicities of
weights in the irreducible representations of the classical Lie
algebras. In the general case, these polynomials have been defined
by Lusztig \cite{Lu1}. For the root system $A_n$, they coincide with
the Kostka-Foulkes polynomials $K_{\lambda\mu}(q)$ (\cf \cite{Mcd}),
which are the coefficients of the Schur symmetric functions
$s_\lambda(X)$ on the basis of Hall-Littlewood functions $P_\mu(X;q)$.
As recently shown by Kirillov and Reshetikhin \cite{KR}, they are
also the generating functions for the sum of quantum numbers of the
Bethe vectors of certain integrable models in Statistical Mechanics.
Also, the specialization at roots of unity of particular Kostka-Foulkes
polynomials gives the decomposition coefficients of certain plethysms
\cite{LLT4}.

A combinatorial interpretation of the Kostka-Foulkes polynomials
has been given in \cite{LS2} (see also \cite{Sc},\cite{Mcd}), where
they were identified as the Poincar\'e polynomials of natural
subsets of the {\it plactic monoid} regarded as a ranked poset.

The plactic monoid is the multiplicative structure on the set of Young tableaux
which reflects the Robinson-Schensted correspondence.

As observed by Date-Jimbo-Miwa \cite{DJM}, the tensor products of the
Gelfand-Tsetlin bases of the quantum enveloping algebra $U_q(\sl_{n+1})$
at $q=0$ are also described by the Robinson-Schensted correspondence.

This was the starting point of Kashiwara's theory of {\it crystal bases},
which are canonical bases of the integrable modules of the algebras
$U_q(\g)$, $\g$ being any symmetrizable Kac-Moody algebra.

The action of the generators of $U_q(\g)$ on the crystal basis at
$q=0$ is described by a combinatorial object, called the {\it crystal
graph} \cite{Ka1}. This reduces a large part of the representation theory
of $U_q(\g)$ to combinatorial questions. For example, one can
read the resolution of a tensor product into its irreducible components
by looking at the connected components of its crystal graph.

The main theorem of the present paper gives a description of the
$q$-multiplicities $K_{\lambda\mu}(q)$, for the root system $A_n$,
in terms of the geometry of the crystal graph associated to
the irreducible $U_q(\sl_{n+1})$-module
$V_\lambda$. The essential tool is an operation of the Weyl group
on the crystal graph, which for $A_n$  coincides
with the canonical action
of the symmetric group on the plactic monoid defined in \cite{LS1}.
The fixed points of this action play a particular r\^ole, and lead
to the definition of multivariate polynomials, which refine the
generating functions of Kostant's generalized exponents for
$SL(n+1,\C)$.

The main results have been announced in \cite{LLT4}.

We would like to thank Andrei Zelevinsky for many stimulating discussions.

\section{Crystal graphs}

Recall that the quantum enveloping algebra $U_q(\sl_{n+1})$ is the
$\Q(q)$-algebra generated by elements $e_i,f_i,t_i,t_i^{-1}$,
$1\le i \le n$, subject
to the following relations \cite{Ji}\cite{Dr}
\begin{equation}
[t_i,t_j] = 0 \ , \qquad t_i t_i^{-1}=t_i^{-1}t_i = 1
\end{equation}
\begin{equation}
[t_i,e_j] = [t_i,f_j] = 0 \ \ {\rm for} \ |i-j|>1
\end{equation}
\begin{equation}
t_je_i=q^{-1}e_it_j \ , \ \ t_jf_i=qf_it_j \ \ {\rm for}\ |i-j|=1
\end{equation}
\begin{equation}
t_ie_i = q^2 e_it_i \ , \ \ \  t_if_i = q^{-2}f_it_i
\end{equation}
\begin{equation}
[e_i,f_j] = \delta_{ij} {t_i-t_i^{-1}\over q-q^{-1}}
\end{equation}
\begin{equation}
[e_i,e_j] = [f_i,f_j] = 0 \ \ {\rm for} \ |i-j|>1
\end{equation}
\begin{equation}
e_je_i^2 -(q+q^{-1})e_ie_je_i +e_i^2e_j =
f_jf_i^2 -(q+q^{-1})f_if_jf_i +f_i^2f_j =0 \ {\rm for} \ |i-j|=1  \ .
\end{equation}
For a $U_q(\sl_{n+1})$-module $M$ and $\lambda\in\Z^{n+1}$, set
$M_\lambda=\{u\in M\ |\ t_i u=q^{\lambda_i-\lambda_{i+1}}u\}$. Elements
of $M_\lambda$ are called weight vectors (of weight $\lambda$), and
a weight  vector is said to be primitive if it is annihilated by
the $e_i$'s. A highest weight $U_q(\sl_{n+1})$-module is a module
$M$ containing a primitive vector $v$ such that $M=U_q(\sl_{n+1})\,v$.
We denote by $V_\lambda$ the irreducible highest weight module
with highest weight $\lambda$.

The module $M$ is said to be integrable if each $M_\lambda$ is
finite-dimensional,
$M=\bigoplus_\lambda M_\lambda$, and for any $i$, $M$ is a direct sum
of finite dimensional representations of the subalgebra isomorphic to
$U_q(\sl_2)$ generated
by $e_i,f_i,t_i$ and $t_i^{-1}$.
By the representation theory of $U_q(\sl_2)$, for any integrable
$U_q(\sl_{n+1})$-module $M$, one has the decomposition
\begin{equation}
M=\bigoplus_\lambda \bigoplus_{0\le k\le \lambda_i-\lambda_{i+1}}
f_i^{(k)} (\ker e_i \cap M_\lambda)
\end{equation}
where $f_i^{(k)}= {q-q^{-1}\over q^k-q^{-k}}\, f_i^k$. Then,
Kashiwara \cite{Ka1} defines endomorphisms $\tilde e_i$, $\tilde f_i \in \End
M$
by
\begin{equation}
\tilde f_i (f_i^{(k)}u)=f_i^{(k+1)} u \ \ {\rm and}\ \
\tilde e_i (f_i^{(k)}u)=f_i^{(k-1)} u
\end{equation}
for $u\in \ker (e_i)\cap M_\lambda$ and $0\le k\le \lambda_i-\lambda_{i+1}$
(by convention $f_i^{(k)} = 0$ for $k<0$).

Let ${\cal A}$ be the subalgebra of $\Q(q)$ formed by the rational
functions without pole at $q=0$. Kashiwara  introduces the
${\cal A}$-lattice in $V_\lambda$
\begin{equation}
L(\lambda)= \sum_{i_1,i_2,\ldots,i_r} {\cal A}
\, \tilde f_{i_1}\tilde f_{i_2}\cdots \tilde f_{i_r}\, u_\lambda
\end{equation}
where $u_\lambda$ is a highest weight vector of $V_\lambda$, and
shows that the set
\begin{equation}
B(\lambda) = \{
\tilde f_{i_1}\tilde f_{i_2}\cdots \tilde f_{i_r}\, u_\lambda \mod
qL(\lambda)\, |\,
1\le i_1,\ldots,i_r\le n\}\backslash\{0\}
\end{equation}
is a basis of the $\Q$-vector space $L(\lambda)/qL(\lambda)$.
The pair $(L(\lambda),B(\lambda))$ is called a (lower) {\it crystal basis}
of $V_\lambda$.
One has
\begin{equation}
\tilde e_i L(\lambda)\subset L(\lambda)\,,\ \
\tilde f_i L(\lambda)\subset L(\lambda)\,,
\end{equation}
so that $\tilde e_i$ and $\tilde f_i$ induce operators in
$L(\lambda)/qL(\lambda)$
still denoted by $\tilde e_i$, $\tilde f_i$. Their action on $B(\lambda)$ is
strikingly simple, namely,
\begin{equation}
\tilde e_i B(\lambda)\subset B(\lambda)\cup \{0\}\,,\ \
\tilde f_i  B(\lambda)\subset B(\lambda)\cup \{0\}\,,
\end{equation}
and for $u,v\in B(\lambda)$,
\begin{equation}
\tilde f_iu = v \ \Longleftrightarrow \ \tilde e_iv = u\,.
\end{equation}
The {\it crystal graph} $\Gamma_\lambda$ associated to $V_\lambda$ is
the coloured oriented graph whose vertices are the elements of $B(\lambda)$,
and whose arrows of colour $i$ describe the action of $\tilde f_i$:
$$
u \stackrel{i}{\longrightarrow} v \ \Longleftrightarrow \ \tilde f_iu = v\,.
$$

More generally,  there exist crystal bases and crystal graphs for any
integrable
module $M$ with highest weights. These objects are compatible with
tensor products, in the sense
that if $(L_1,B_1)$ and $(L_2,B_2)$ are crystal bases for the modules
$M_1$ and $M_2$ then $(L_1\otimes L_2,B_1\otimes B_2)$ is a crystal
basis of $M_1\otimes M_2$. Moreover, the decomposition of $M_1\otimes M_2$
into irreducible representations is given by the decomposition
of its crystal graph into  connected components,
two submodules  being isomorphic if and only if
the corresponding  graphs are isomorphic (as coloured graphs).
In particular, applying this process to the $r$-th tensor power
of the fundamental representation $V$ whose crystal graph is
$$
1 \stackrel{1}{\longrightarrow}
2 \stackrel{2}{\longrightarrow}
\cdots
n-1 \stackrel{n-1}{\longrightarrow}
n \stackrel{n}{\longrightarrow}
n+1
$$
(where for brevity the basis vector $v_i$ is denoted by
its index $i$), one obtains by identifying the isomorphic
connected components an equivalence relation on the set of words labeling
the vertices of the crystal graph of $V^{\otimes r}$. This relation
is  the plactic equivalence, described in the
next section.

\section{The plactic monoid} \label{PLACTIC}

Let $A=\{a_1<a_2< \ldots <a_{n+1}\}$ denote a totally ordered
alphabet of noncommutative indeterminates,
and consider the free monoid $A^*$ generated by $A$.
The Robinson-Schensted correspondence
associates to a word $w \in A^*$ a pair $(P(w),Q(w))$ of Young tableaux of the
same
shape. It was shown by Knuth \cite{Kn} that the equivalence on $A^*$ defined by
$$ w \equiv w'\ \Longleftrightarrow \ P(w) = P(w') $$
is generated by the relations
$$zxy \equiv xzy \ , \ \ \ yxz \equiv yzx  \ , \ \ \
yxx \equiv xyx \ , \ \ \ yxy \equiv yyx  \ , $$
for any $x{<}y{<}z$ in $A$. The quotient set $A^*/{\equiv}$, which by
definition is in
one-to-one correspondence with the set of Young tableaux on $A$, is therefore
endowed
with a multiplicative structure reflecting the Robinson-Schensted construction.
The monoid $A^*/{\equiv}$ is called the {\it plactic monoid}.
We denote by $\tab (\lambda,\mu)$ the set of tableaux with shape
$\lambda$ and weight $\mu$, \ie such that the number of occurences
of the letter $a_i$ is $\mu_i$.

As can be deduced from
the detailed description of crystal graphs given in
\cite{KN}, two vertices $v_{i_1}\otimes v_{i_2}\otimes\cdots\otimes v_{i_r}$
and $v_{j_1}\otimes v_{j_2}\otimes\cdots\otimes v_{j_r}$ of the
crystal graph of $V^{\otimes r}$ are identified (for the equivalence
mentioned at the end of the preceding section) if and only if the
words $a_{i_1}a_{i_2}\cdots a_{i_r}$ and $a_{j_1}a_{j_2}\cdots a_{j_r}$
are plactically congruent. In fact, the Knuth relations amount
to the identification of the two copies of $V_{(2,1)}$ lying in $V^{\otimes
3}$.

The algebraic properties of the plactic monoid have been investigated in
\cite{LS1}.

Although this will not be used in the sequel, it is worthwhile
to mention at this point that the above approach to the
plactic monoid allows to define a similar object for any finite
dimensional complex simple Lie algebra. For example, the plactic
monoid of type $C_n$ is the quotient of the free monoid on the
ordered alphabet
$$\bar A + A=
\{\bar a_n < \bar a_{n-1} <\ldots <\bar a_1 < a_1 <\ldots < a_n\}$$
by the relations (\ref{SP1}) to (\ref{SP4})

\begin{equation}\label{SP1}
{\rm for}\ \ x<y<z \ \ {\rm and}\ \ x\not=\bar z\ ,\qquad
\left\{ \matrix{
yxz & \equiv & yzx \cr
xzy & \equiv & zxy \cr}\right.
\end{equation}
\begin{equation}\label{SP2}
{\rm for}\ \ x<y \ \ {\rm and}\ \ x\not=\bar y\ ,\qquad
\left\{ \matrix{
yxx & \equiv & xyx \cr
yyx & \equiv & yxy \cr}\right.
\end{equation}
\begin{equation}\label{SP3}
{\rm for}\ \ i\le n-1\ \ {\rm and}\ \  \bar a_i \le x \le a_i  \ ,\qquad
\left\{\matrix{
a_i\bar a_i x & \equiv & \bar a_{i+1} a_{i+1} x \cr
x a_i \bar a_i & \equiv & x \bar a_{i+1} a_{i+1} \cr}\right.
\end{equation}
For $w=x_1x_2\cdots x_k$ with $x_1>x_2>\ldots >x_k$, set
$$
E_w =\{ (x_p,x_q)\ |\ x_p=a_i \, ,\
x_q = \bar a_i\, ,\ q-p < k+i-n \} \ .
$$
Then,
\begin{equation}\label{SP4}
w \equiv \hat w
\end{equation}
where $\hat w$ is the word obtained from $w$ by erasing all pairs
$(x_p,x_q)\in E_w$.
These relations are obtained by identifying the isomorphic connected
components of the crystal graph of $V^{\otimes 3}$ where $V$ is
the vector representation of $U_q(\sp_{2n})$. Distinguished representatives
of the equivalence classes are provided by the tableaux of \cite{KN},
up to a simple reindexation due to the use of a different ordering.
One can deduce from this a  Robinson-Schensted type algorithm, sending a word
$w$ on $\bar A +A$ onto a pair $(P(w),Q(w))$, where $P(w)$ is
the above mentioned representative of the class of $w$, and $Q(w)$
is an oscillating tableau, describing the intermediate shapes
arising at the succesive stages of the algorithm. This bijection
is different from the one given by Berele \cite{Ber}.

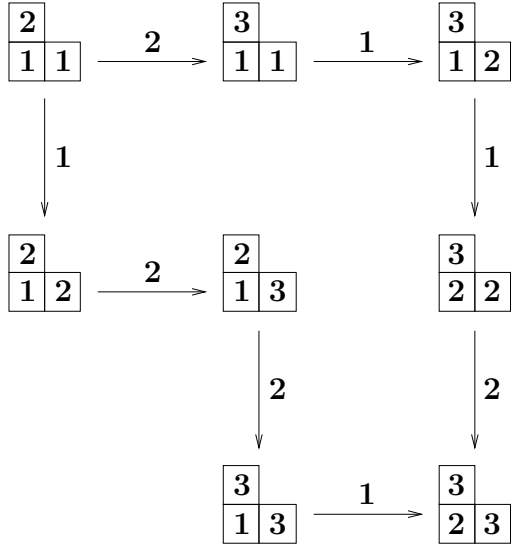
\begin{figure}[t]
\begin{center}
\setlength{\unitlength}{0.0085in}
\begingroup\makeatletter\ifx\SetFigFont\undefined
\def\x#1#2#3#4#5#6#7\relax{\def\x{#1#2#3#4#5#6}}%
\expandafter\x\fmtname xxxxxx\relax \def\y{splain}%
\ifx\x\y   
\gdef\SetFigFont#1#2#3{%
  \ifnum #1<17\tiny\else \ifnum #1<20\small\else
  \ifnum #1<24\normalsize\else \ifnum #1<29\large\else
  \ifnum #1<34\Large\else \ifnum #1<41\LARGE\else
     \huge\fi\fi\fi\fi\fi\fi
  \csname #3\endcsname}%
\else
\gdef\SetFigFont#1#2#3{\begingroup
  \count@#1\relax \ifnum 25<\count@\count@25\fi
  \def\x{\endgroup\@setsize\SetFigFont{#2pt}}%
  \expandafter\x
    \csname \romannumeral\the\count@ pt\expandafter\endcsname
    \csname @\romannumeral\the\count@ pt\endcsname
  \csname #3\endcsname}%
\fi
\fi\endgroup
\begin{picture}(309,349)(0,-10)
\path(22,309)(22,333)(0,333)
	(0,309)(22,309)
\path(22,285)(22,309)(0,309)
	(0,285)(22,285)
\path(44,285)(44,309)(22,309)
	(22,285)(44,285)
\path(22,167)(22,190)(0,190)
	(0,167)(22,167)
\path(44,143)(44,167)(22,167)
	(22,143)(44,143)
\path(22,143)(22,167)(0,167)
	(0,143)(22,143)
\path(154,309)(154,333)(132,333)
	(132,309)(154,309)
\path(154,285)(154,309)(132,309)
	(132,285)(154,285)
\path(177,285)(177,309)(154,309)
	(154,285)(177,285)
\path(287,309)(287,333)(265,333)
	(265,309)(287,309)
\path(287,285)(287,309)(265,309)
	(265,285)(287,285)
\path(309,285)(309,309)(287,309)
	(287,285)(309,285)
\path(154,167)(154,190)(132,190)
	(132,167)(154,167)
\path(154,143)(154,167)(132,167)
	(132,143)(154,143)
\path(177,143)(177,167)(154,167)
	(154,143)(177,143)
\path(287,167)(287,190)(265,190)
	(265,167)(287,167)
\path(287,143)(287,167)(265,167)
	(265,143)(287,143)
\path(309,143)(309,167)(287,167)
	(287,143)(309,143)
\path(154,24)(154,48)(132,48)
	(132,24)(154,24)
\path(154,0)(154,24)(132,24)
	(132,0)(154,0)
\path(177,0)(177,24)(154,24)
	(154,0)(177,0)
\path(287,24)(287,48)(265,48)
	(265,24)(287,24)
\path(287,0)(287,24)(265,24)
	(265,0)(287,0)
\path(309,0)(309,24)(287,24)
	(287,0)(309,0)
\path(55,297)(121,297)
\path(113.000,295.000)(121.000,297.000)(113.000,299.000)
\path(188,297)(254,297)
\path(246.000,295.000)(254.000,297.000)(246.000,299.000)
\path(287,274)(287,202)
\path(285.000,210.000)(287.000,202.000)(289.000,210.000)
\path(287,131)(287,59)
\path(285.000,67.000)(287.000,59.000)(289.000,67.000)
\path(188,18)(254,18)
\path(246.000,16.000)(254.000,18.000)(246.000,20.000)
\path(154,131)(154,59)
\path(152.000,67.000)(154.000,59.000)(156.000,67.000)
\path(55,155)(121,155)
\path(113.000,153.000)(121.000,155.000)(113.000,157.000)
\path(22,274)(22,202)
\path(20.000,210.000)(22.000,202.000)(24.000,210.000)
\put(6,315){\makebox(0,0)[lb]{\smash{{{\SetFigFont{12}{14.4}{bf}2}}}}}
\put(6,291){\makebox(0,0)[lb]{\smash{{{\SetFigFont{12}{14.4}{bf}1}}}}}
\put(28,291){\makebox(0,0)[lb]{\smash{{{\SetFigFont{12}{14.4}{bf}1}}}}}
\put(138,315){\makebox(0,0)[lb]{\smash{{{\SetFigFont{12}{14.4}{bf}3}}}}}
\put(138,291){\makebox(0,0)[lb]{\smash{{{\SetFigFont{12}{14.4}{bf}1}}}}}
\put(160,291){\makebox(0,0)[lb]{\smash{{{\SetFigFont{12}{14.4}{bf}1}}}}}
\put(270,315){\makebox(0,0)[lb]{\smash{{{\SetFigFont{12}{14.4}{bf}3}}}}}
\put(270,291){\makebox(0,0)[lb]{\smash{{{\SetFigFont{12}{14.4}{bf}1}}}}}
\put(292,291){\makebox(0,0)[lb]{\smash{{{\SetFigFont{12}{14.4}{bf}2}}}}}
\put(6,172){\makebox(0,0)[lb]{\smash{{{\SetFigFont{12}{14.4}{bf}2}}}}}
\put(6,149){\makebox(0,0)[lb]{\smash{{{\SetFigFont{12}{14.4}{bf}1}}}}}
\put(28,149){\makebox(0,0)[lb]{\smash{{{\SetFigFont{12}{14.4}{bf}2}}}}}
\put(138,149){\makebox(0,0)[lb]{\smash{{{\SetFigFont{12}{14.4}{bf}1}}}}}
\put(160,149){\makebox(0,0)[lb]{\smash{{{\SetFigFont{12}{14.4}{bf}3}}}}}
\put(270,172){\makebox(0,0)[lb]{\smash{{{\SetFigFont{12}{14.4}{bf}3}}}}}
\put(270,149){\makebox(0,0)[lb]{\smash{{{\SetFigFont{12}{14.4}{bf}2}}}}}
\put(292,149){\makebox(0,0)[lb]{\smash{{{\SetFigFont{12}{14.4}{bf}2}}}}}
\put(270,30){\makebox(0,0)[lb]{\smash{{{\SetFigFont{12}{14.4}{bf}3}}}}}
\put(270,6){\makebox(0,0)[lb]{\smash{{{\SetFigFont{12}{14.4}{bf}2}}}}}
\put(292,6){\makebox(0,0)[lb]{\smash{{{\SetFigFont{12}{14.4}{bf}3}}}}}
\put(138,30){\makebox(0,0)[lb]{\smash{{{\SetFigFont{12}{14.4}{bf}3}}}}}
\put(138,6){\makebox(0,0)[lb]{\smash{{{\SetFigFont{12}{14.4}{bf}1}}}}}
\put(160,6){\makebox(0,0)[lb]{\smash{{{\SetFigFont{12}{14.4}{bf}3}}}}}
\put(83,303){\makebox(0,0)[lb]{\smash{{{\SetFigFont{12}{14.4}{bf}2}}}}}
\put(215,303){\makebox(0,0)[lb]{\smash{{{\SetFigFont{12}{14.4}{bf}1}}}}}
\put(292,232){\makebox(0,0)[lb]{\smash{{{\SetFigFont{12}{14.4}{bf}1}}}}}
\put(292,89){\makebox(0,0)[lb]{\smash{{{\SetFigFont{12}{14.4}{bf}2}}}}}
\put(215,24){\makebox(0,0)[lb]{\smash{{{\SetFigFont{12}{14.4}{bf}1}}}}}
\put(160,89){\makebox(0,0)[lb]{\smash{{{\SetFigFont{12}{14.4}{bf}2}}}}}
\put(83,161){\makebox(0,0)[lb]{\smash{{{\SetFigFont{12}{14.4}{bf}2}}}}}
\put(28,232){\makebox(0,0)[lb]{\smash{{{\SetFigFont{12}{14.4}{bf}1}}}}}
\put(138,172){\makebox(0,0)[lb]{\smash{{{\SetFigFont{12}{14.4}{bf}2}}}}}
\end{picture}
\caption{Crystal graph of the $U_q(\sl_3)$-module $V_{(2,\,1)}$}
\end{center}
\end{figure}

Now returning to the $A_n$ case, we introduce following \cite{LS1,LS3,LS4}
linear operators $\epsilon_i,\, \phi_i,\, \sigma_i,\ i{=}1,\,\ldots ,\,n$
acting on the free associative algebra $\Q\<A\>$ in the following way.
Consider first the case of a two-letter alphabet $A=\{a_i,a_{i+1}\}$. Let
$w = x_1\cdots x_m$ be
a word on $A$. Bracket every factor $a_{i+1}\, a_i$ of $w$. The letters which
are
not bracketed constitute a subword $w_1$ of $w$. Then, bracket every factor
$a_{i+1}\, a_i$ of $w_1$. There remains a subword $w_2$. Continue this
procedure
until it stops, giving a word $w_k=x_{j_1}\cdots x_{j_{r+s}}$ of the form
$w_k=a_i^r\, a_{i+1}^s$ for some integers $r,\,s$.
The image of $w_k$ under $\epsilon_i,\, \phi_i$ or $\sigma_i$ is given by
$$
\epsilon_i(a_i^r a_{i+1}^s)
=\left\{ \matrix{ a_i^{r+1} a_{i+1}^{s-1} & (s\ge 1) \cr
                                                  0               &
(s=0)}\right. ,\ \
 \phi_i(a_i^r a_{i+1}^s)=
 \left\{ \matrix{ a_i^{r-1} a_{i+1}^{s+1} & (r\ge 1) \cr
          0               & (r=0)}\right. ,\ \
$$
$$
 \sigma_i(a_i^r a_{i+1}^s)= a_i^s a_{i+1}^r\,.
$$
The image of the initial word $w$ is then $w' = y_1\cdots y_m$ where
$y_k=x_k$ for $k \not \in \{j_1\,,\ldots ,\,j_{r+s}\}$.
For instance, if
$$w=(a_2a_1)a_1a_1a_2(a_2a_1)a_1a_1a_1a_2\,,$$
we shall have
$$w_1=.\ .\ a_1a_1(a_2\ . \ . \ a_1)a_1a_1a_2\,,$$
$$w_2=.\ .\ a_1a_1\ .\ .\ .\ .\ a_1a_1a_2\,.$$
Thus, $$\epsilon_1(w)=a_2a_1{\bf a_1a_1}a_2a_2a_1a_1{\bf a_1a_1a_1}\,,$$
$$\phi_1(w)=a_2a_1{\bf a_1a_1}a_2a_2a_1a_1{\bf a_1a_2a_2}\,,$$
$$\sigma_1(w)=a_2a_1{\bf a_1a_2}a_2a_2a_1a_1{\bf a_2a_2a_2}\,,$$
where the letters printed in bold type are those of the image of the subword
$w_2$. Finally, in the general case, the action of the operators
$\epsilon_i,\, \phi_i,\, \sigma_i$ on $w$ is defined by the previous
rules applied to the subword consisting of the letters $a_i,\,a_{i+1}$,
the remaining letters being unchanged.
We also define an involution denoted by $\Omega_2$. It is the
{\it anti}-automorphism of the algebra $\Q\<A\>$ such that
$ \Omega_2 (a_i) = a_{n-i+2}$.

All these operators are compatible with the plactic congruence ${\equiv}$
and therefore induce operations on tableaux that we shall also
denote by $\epsilon_i,\, \phi_i,\, \sigma_i$ and $\Omega_2$.

Let $\lambda$ be a partition of length $\le n+1$, and let $y_{\lambda}$
be the corresponding Yamanouchi tableau, \ie the only tableau of shape
and weight $\lambda$. It is known \cite{LS4} that the products
$\phi_{i_1}\cdots \phi_{i_r}$ applied to $y_{\lambda}$
generate the set $\tab(\lambda\, ,\,\cdot)$ of all tableaux of shape $\lambda$
over $A$. Define
a coloured graph $G_\lambda$ on $\tab(\lambda\, ,\,\cdot)$
whose edges of colour $i$ describe
the action of~$\phi_i$:
$$
t  \stackrel{i}{\longrightarrow} t' \ \  \Longleftrightarrow \ \  \phi_i(t) =
t'\,.
$$
Comparing this definition with \cite{KN}, one sees that $G_\lambda$ is
none other than the crystal graph $\Gamma_\lambda$ of $V_\lambda$.
In other words, the operators $\epsilon_i$, $\phi_i$ coincide with
the endomorphisms $\tilde e_i$, $\tilde f_i$ of $V^{\otimes r}$
``at $q=0$", in the identification
$a_{i_1}\cdots a_{i_r} \longleftrightarrow v_{i_1}\otimes \cdots \otimes
v_{i_r}$.

It follows from the definition of $\Omega_2$ that
$$
t  \stackrel{i}{\longrightarrow} t' \ \  \Longleftrightarrow \ \
\Omega_2(t')  \stackrel{n+1-i}{\longrightarrow} \Omega_2(t) \,.
$$
Hence, the crystal graph has a symmetry of order 2 obtained by
reversing its arrows and changing their colour from $i$ to $n+1-i$.

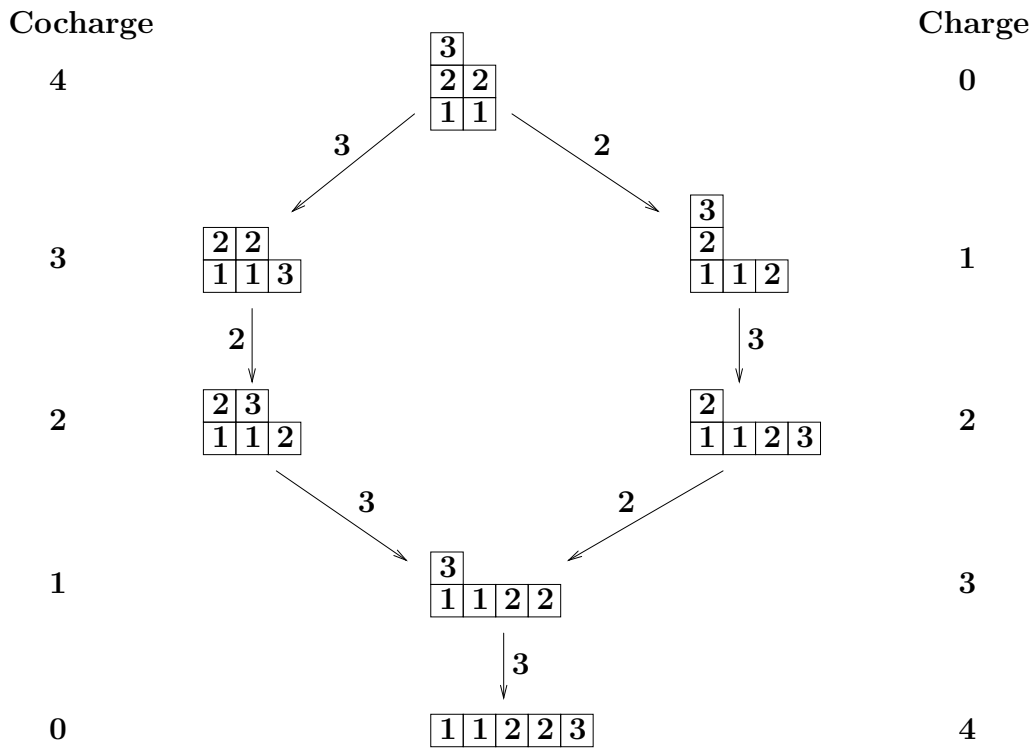
\begin{figure}[t]
\begin{center}
\setlength{\unitlength}{0.0085 in}
\begingroup\makeatletter\ifx\SetFigFont\undefined
\def\x#1#2#3#4#5#6#7\relax{\def\x{#1#2#3#4#5#6}}%
\expandafter\x\fmtname xxxxxx\relax \def\y{splain}%
\ifx\x\y   
\gdef\SetFigFont#1#2#3{%
  \ifnum #1<17\tiny\else \ifnum #1<20\small\else
  \ifnum #1<24\normalsize\else \ifnum #1<29\large\else
  \ifnum #1<34\Large\else \ifnum #1<41\LARGE\else
     \huge\fi\fi\fi\fi\fi\fi
  \csname #3\endcsname}%
\else
\gdef\SetFigFont#1#2#3{\begingroup
  \count@#1\relax \ifnum 25<\count@\count@25\fi
  \def\x{\endgroup\@setsize\SetFigFont{#2pt}}%
  \expandafter\x
    \csname \romannumeral\the\count@ pt\expandafter\endcsname
    \csname @\romannumeral\the\count@ pt\endcsname
  \csname #3\endcsname}%
\fi
\fi\endgroup
\begin{picture}(604,472)(0,-10)
\path(440,320)(440,340)(420,340)
	(420,320)(440,320)
\path(440,300)(440,320)(420,320)
	(420,300)(440,300)
\path(440,280)(440,300)(420,300)
	(420,280)(440,280)
\path(460,280)(460,300)(440,300)
	(440,280)(460,280)
\path(480,280)(480,300)(460,300)
	(460,280)(480,280)
\path(440,200)(440,220)(420,220)
	(420,200)(440,200)
\path(440,180)(440,200)(420,200)
	(420,180)(440,180)
\path(460,180)(460,200)(440,200)
	(440,180)(460,180)
\path(480,180)(480,200)(460,200)
	(460,180)(480,180)
\path(500,180)(500,200)(480,200)
	(480,180)(500,180)
\path(280,420)(280,440)(260,440)
	(260,420)(280,420)
\path(280,400)(280,420)(260,420)
	(260,400)(280,400)
\path(300,400)(300,420)(280,420)
	(280,400)(300,400)
\path(280,380)(280,400)(260,400)
	(260,380)(280,380)
\path(300,380)(300,400)(280,400)
	(280,380)(300,380)
\path(140,300)(140,320)(120,320)
	(120,300)(140,300)
\path(160,300)(160,320)(140,320)
	(140,300)(160,300)
\path(140,280)(140,300)(120,300)
	(120,280)(140,280)
\path(160,280)(160,300)(140,300)
	(140,280)(160,280)
\path(180,280)(180,300)(160,300)
	(160,280)(180,280)
\path(140,200)(140,220)(120,220)
	(120,200)(140,200)
\path(160,200)(160,220)(140,220)
	(140,200)(160,200)
\path(140,180)(140,200)(120,200)
	(120,180)(140,180)
\path(160,180)(160,200)(140,200)
	(140,180)(160,180)
\path(180,180)(180,200)(160,200)
	(160,180)(180,180)
\path(280,100)(280,120)(260,120)
	(260,100)(280,100)
\path(280,80)(280,100)(260,100)
	(260,80)(280,80)
\path(300,80)(300,100)(280,100)
	(280,80)(300,80)
\path(320,80)(320,100)(300,100)
	(300,80)(320,80)
\path(340,80)(340,100)(320,100)
	(320,80)(340,80)
\path(280,0)(280,20)(260,20)
	(260,0)(280,0)
\path(300,0)(300,20)(280,20)
	(280,0)(300,0)
\path(320,0)(320,20)(300,20)
	(300,0)(320,0)
\path(340,0)(340,20)(320,20)
	(320,0)(340,0)
\path(360,0)(360,20)(340,20)
	(340,0)(360,0)
\path(250,390)(175,330)
\path(179.998,336.559)(175.000,330.000)(182.496,333.436)
\path(150,270)(150,225)
\path(148.000,233.000)(150.000,225.000)(152.000,233.000)
\path(165,170)(245,115)
\path(237.275,117.884)(245.000,115.000)(239.541,121.180)
\path(310,390)(400,330)
\path(392.234,332.774)(400.000,330.000)(394.453,336.102)
\path(450,270)(450,225)
\path(448.000,233.000)(450.000,225.000)(452.000,233.000)
\path(440,170)(345,115)
\path(350.921,120.739)(345.000,115.000)(352.925,117.277)
\path(305,70)(305,30)
\path(303.000,38.000)(305.000,30.000)(307.000,38.000)
\put(265,425){\makebox(0,0)[lb]{\smash{{{\SetFigFont{12}{14.4}{bf}3}}}}}
\put(265,405){\makebox(0,0)[lb]{\smash{{{\SetFigFont{12}{14.4}{bf}2}}}}}
\put(285,405){\makebox(0,0)[lb]{\smash{{{\SetFigFont{12}{14.4}{bf}2}}}}}
\put(265,385){\makebox(0,0)[lb]{\smash{{{\SetFigFont{12}{14.4}{bf}1}}}}}
\put(285,385){\makebox(0,0)[lb]{\smash{{{\SetFigFont{12}{14.4}{bf}1}}}}}
\put(125,305){\makebox(0,0)[lb]{\smash{{{\SetFigFont{12}{14.4}{bf}2}}}}}
\put(145,305){\makebox(0,0)[lb]{\smash{{{\SetFigFont{12}{14.4}{bf}2}}}}}
\put(125,285){\makebox(0,0)[lb]{\smash{{{\SetFigFont{12}{14.4}{bf}1}}}}}
\put(145,285){\makebox(0,0)[lb]{\smash{{{\SetFigFont{12}{14.4}{bf}1}}}}}
\put(165,285){\makebox(0,0)[lb]{\smash{{{\SetFigFont{12}{14.4}{bf}3}}}}}
\put(125,205){\makebox(0,0)[lb]{\smash{{{\SetFigFont{12}{14.4}{bf}2}}}}}
\put(145,205){\makebox(0,0)[lb]{\smash{{{\SetFigFont{12}{14.4}{bf}3}}}}}
\put(125,185){\makebox(0,0)[lb]{\smash{{{\SetFigFont{12}{14.4}{bf}1}}}}}
\put(145,185){\makebox(0,0)[lb]{\smash{{{\SetFigFont{12}{14.4}{bf}1}}}}}
\put(165,185){\makebox(0,0)[lb]{\smash{{{\SetFigFont{12}{14.4}{bf}2}}}}}
\put(265,105){\makebox(0,0)[lb]{\smash{{{\SetFigFont{12}{14.4}{bf}3}}}}}
\put(265,85){\makebox(0,0)[lb]{\smash{{{\SetFigFont{12}{14.4}{bf}1}}}}}
\put(285,85){\makebox(0,0)[lb]{\smash{{{\SetFigFont{12}{14.4}{bf}1}}}}}
\put(305,85){\makebox(0,0)[lb]{\smash{{{\SetFigFont{12}{14.4}{bf}2}}}}}
\put(325,85){\makebox(0,0)[lb]{\smash{{{\SetFigFont{12}{14.4}{bf}2}}}}}
\put(265,5){\makebox(0,0)[lb]{\smash{{{\SetFigFont{12}{14.4}{bf}1}}}}}
\put(285,5){\makebox(0,0)[lb]{\smash{{{\SetFigFont{12}{14.4}{bf}1}}}}}
\put(305,5){\makebox(0,0)[lb]{\smash{{{\SetFigFont{12}{14.4}{bf}2}}}}}
\put(325,5){\makebox(0,0)[lb]{\smash{{{\SetFigFont{12}{14.4}{bf}2}}}}}
\put(345,5){\makebox(0,0)[lb]{\smash{{{\SetFigFont{12}{14.4}{bf}3}}}}}
\put(425,325){\makebox(0,0)[lb]{\smash{{{\SetFigFont{12}{14.4}{bf}3}}}}}
\put(425,305){\makebox(0,0)[lb]{\smash{{{\SetFigFont{12}{14.4}{bf}2}}}}}
\put(425,285){\makebox(0,0)[lb]{\smash{{{\SetFigFont{12}{14.4}{bf}1}}}}}
\put(445,285){\makebox(0,0)[lb]{\smash{{{\SetFigFont{12}{14.4}{bf}1}}}}}
\put(465,285){\makebox(0,0)[lb]{\smash{{{\SetFigFont{12}{14.4}{bf}2}}}}}
\put(425,205){\makebox(0,0)[lb]{\smash{{{\SetFigFont{12}{14.4}{bf}2}}}}}
\put(425,185){\makebox(0,0)[lb]{\smash{{{\SetFigFont{12}{14.4}{bf}1}}}}}
\put(445,185){\makebox(0,0)[lb]{\smash{{{\SetFigFont{12}{14.4}{bf}1}}}}}
\put(465,185){\makebox(0,0)[lb]{\smash{{{\SetFigFont{12}{14.4}{bf}2}}}}}
\put(485,185){\makebox(0,0)[lb]{\smash{{{\SetFigFont{12}{14.4}{bf}3}}}}}
\put(200,365){\makebox(0,0)[lb]{\smash{{{\SetFigFont{12}{14.4}{bf}3}}}}}
\put(360,365){\makebox(0,0)[lb]{\smash{{{\SetFigFont{12}{14.4}{bf}2}}}}}
\put(455,245){\makebox(0,0)[lb]{\smash{{{\SetFigFont{12}{14.4}{bf}3}}}}}
\put(375,145){\makebox(0,0)[lb]{\smash{{{\SetFigFont{12}{14.4}{bf}2}}}}}
\put(215,145){\makebox(0,0)[lb]{\smash{{{\SetFigFont{12}{14.4}{bf}3}}}}}
\put(310,45){\makebox(0,0)[lb]{\smash{{{\SetFigFont{12}{14.4}{bf}3}}}}}
\put(0,440){\makebox(0,0)[lb]{\smash{{{\SetFigFont{12}{14.4}{bf}Cocharge}}}}}
\put(25,405){\makebox(0,0)[lb]{\smash{{{\SetFigFont{12}{14.4}{bf}4}}}}}
\put(25,295){\makebox(0,0)[lb]{\smash{{{\SetFigFont{12}{14.4}{bf}3}}}}}
\put(25,195){\makebox(0,0)[lb]{\smash{{{\SetFigFont{12}{14.4}{bf}2}}}}}
\put(25,95){\makebox(0,0)[lb]{\smash{{{\SetFigFont{12}{14.4}{bf}1}}}}}
\put(25,5){\makebox(0,0)[lb]{\smash{{{\SetFigFont{12}{14.4}{bf}0}}}}}
\put(560,440){\makebox(0,0)[lb]{\smash{{{\SetFigFont{12}{14.4}{bf}Charge}}}}}
\put(585,405){\makebox(0,0)[lb]{\smash{{{\SetFigFont{12}{14.4}{bf}0}}}}}
\put(585,295){\makebox(0,0)[lb]{\smash{{{\SetFigFont{12}{14.4}{bf}1}}}}}
\put(585,195){\makebox(0,0)[lb]{\smash{{{\SetFigFont{12}{14.4}{bf}2}}}}}
\put(585,95){\makebox(0,0)[lb]{\smash{{{\SetFigFont{12}{14.4}{bf}3}}}}}
\put(585,5){\makebox(0,0)[lb]{\smash{{{\SetFigFont{12}{14.4}{bf}4}}}}}
\put(135,245){\makebox(0,0)[lb]{\smash{{{\SetFigFont{12}{14.4}{bf}2}}}}}
\end{picture}
\caption{The cyclage graph $H_{(2,2,1)}$ for $A_2$\label{POSET}}
\end{center}
\end{figure}

We end this section by defining another coloured graph issued from the
plactic monoid. It would be interesting to have also an interpretation of
this graph in the framework of crystal bases.

The operation of conjugation $h \rightarrow g^{-1}hg$ in
a group is replaced in a monoid by the circular permutation
$w=uv \rightarrow vu$. In the plactic monoid, a special kind of conjugation
called {\it cyclage} has been introduced in \cite{LS1}
(see also \cite{La}). It is defined as follows.
Let $t,\,t'$ be two tableaux of weight $\mu$, where
$\mu$ is a partition, {\it i. e.} $\mu_1 \ge \mu_2 \ge \ldots \ge \mu_{n+1} \ge
0$.
For $i\ge 2$, write $t \stackrel{i}{\longrightarrow} t'$
if and only if there exists $u$ in $A^*/{\equiv}$ such that
$t=a_i\,u$ and $t'=u\,a_i$. In this case,  $t'$ is said
to be a {\it cyclage} of $t$.
In this way, the set $\tab(\cdot\, ,\mu)$ of all tableaux of weight $\mu$
is given the structure
of a connected oriented (coloured) graph $H_\mu$, whose transitive closure is a
partial order with
minimal element the row tableau
 $t_{\mu} := a_1^{\mu_1}a_2^{\mu_2}\cdots a_{n+1}^{\mu_{n+1}}$.
The cocharge $co(t)$ of a tableau $t$ of weight $\mu$ is defined as its rank
in the poset $\tab(\cdot\, , \, \mu)$, that is, the number of cyclages needed
to transform $t$ into the row tableau $t_{\mu}$.
The maximal value of the cocharge on $\tab(\cdot\, , \, \mu)$
is $ \|\mu\|=\sum_i (i-1)\mu_i$. The {\it charge} of $t$ is $c(t) := \|\mu\| -
co(t)$.

A tableau $t$ which is not a row tableau admits in general several cyclages.
We call {\it initial cyclage}, and we denote by ${\cal C}(t)$, the cyclage
obtained
by cycling the first letter of the row reading of $t$. In Figure~\ref{POSET},
all arrows except the one labelled $2$ in the top right corner correspond to
initial
cyclages. Deleting from $H_\mu$ all
the arrows corresponding to non initial cyclages, one is left with a tree
$T_\mu$ whose root is the row tableau $t_\mu$.

More generally, the initial cyclage of a tableau $t$ whose weight is not
a partition is defined in exactly the same way. The initial cyclage
commutes with the operators $\sigma_i$ \cite{La}:
$$ \sigma_i({\cal C}(t)) = {\cal C}(\sigma_i(t)) $$
for any tableau $t$ whose shape is not a row. Thus,
the image by $\sigma_i$ of the tree $T_\mu$ is
an isomorphic tree $T_{\sigma_i(\mu)}$, and
it is natural to extend
the definition of the charge to tableaux of arbitrary weight by requiring that
$c(\sigma_i(t)) = c(t)$.

\section{Kostka-Foulkes polynomials}

Let $R$ be the root system of a finite dimensional complex simple Lie
algebra $\g$, $W$ its
Weyl group,
$R^+$ the set of positive roots
and $\rho$ the half-sum of positive roots. One can define a $q$-analogue ${\cal
P}_q$ of
Kostant's partition function by
$$
\prod_{\alpha\in R^+} (1-q e^\alpha)^{-1} = \sum_{\mu\in P}{\cal P}_q(\mu)e^\mu
$$
where $P$ is the weight lattice of $\g$.
Substituting this $q$-analogue to the ordinary partition function
in Kostant's weight multiplicity formula, one defines
$$
K_{\lambda\mu}(q) = \sum_{w\in W} (-1)^{\ell(w)}
{\cal P}_q ((\lambda+\rho)w -\mu-\rho)
$$
which, as shown by Lusztig \cite{Lu1}, turns out to be a polynomial in $q$ with
nonnegative integer coefficients. In the general case, the only known proof
of this property comes from the interpretation of the $K_{\lambda\mu}(q)$
as (renormalized) Kazhdan-Lusztig polynomials for the affine Weyl group.
However, for $R=A_n$, they coincide with the Kostka-Foulkes polynomials,
which are the coefficients of the expansion
\begin{equation}
s_\lambda(X) = \sum_\mu K_{\lambda\mu}(q) P_\mu(X;q)
\end{equation}
of the Schur functions on the basis of Hall-Littlewood functions \cite{Mcd}.
In this case, there exists a combinatorial description
of the $K_{\lambda\mu}(q)$, which implies the positivity of their
coefficients.

\begin{theorem}{\rm \cite{LS2,Sc}}
The Kostka-Foulkes polynomial $K_{\lambda \mu}(q)$ is the generating function
of
the charge on $\tab(\lambda ,\mu)$:
$$K_{\lambda \mu}(q)= \sum_{t\in\tab(\lambda ,\mu)} q^{c(t)} \ .$$
\end{theorem}

\section{Crystal graphs and $q$-multiplicities}

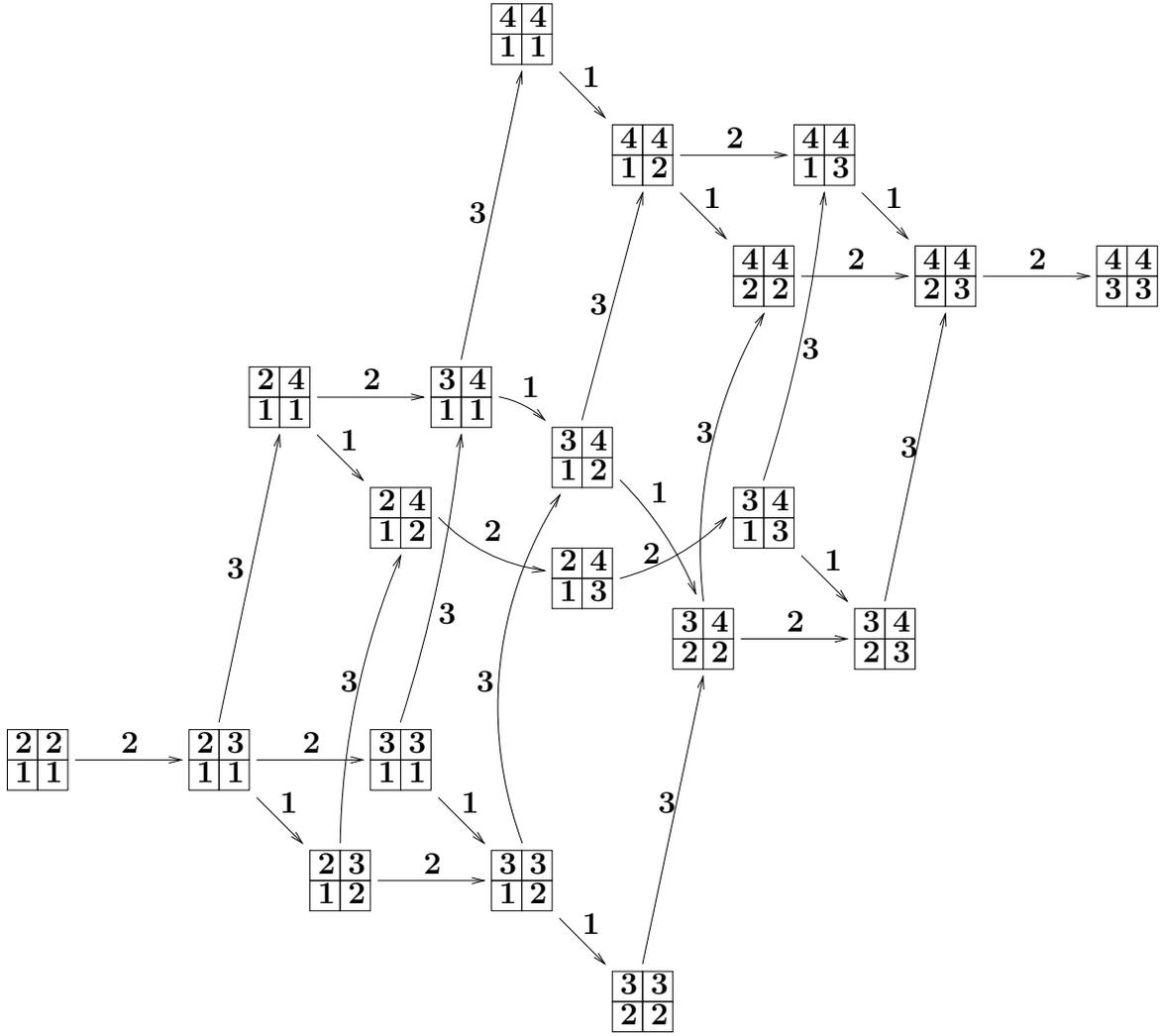
\begin{figure}[t]
\begin{center}
\setlength{\unitlength}{0.008in}
\begingroup\makeatletter\ifx\SetFigFont\undefined
\def\x#1#2#3#4#5#6#7\relax{\def\x{#1#2#3#4#5#6}}%
\expandafter\x\fmtname xxxxxx\relax \def\y{splain}%
\ifx\x\y   
\gdef\SetFigFont#1#2#3{%
  \ifnum #1<17\tiny\else \ifnum #1<20\small\else
  \ifnum #1<24\normalsize\else \ifnum #1<29\large\else
  \ifnum #1<34\Large\else \ifnum #1<41\LARGE\else
     \huge\fi\fi\fi\fi\fi\fi
  \csname #3\endcsname}%
\else
\gdef\SetFigFont#1#2#3{\begingroup
  \count@#1\relax \ifnum 25<\count@\count@25\fi
  \def\x{\endgroup\@setsize\SetFigFont{#2pt}}%
  \expandafter\x
    \csname \romannumeral\the\count@ pt\expandafter\endcsname
    \csname @\romannumeral\the\count@ pt\endcsname
  \csname #3\endcsname}%
\fi
\fi\endgroup
\begin{picture}(760,697)(0,-10)
\path(180,420)(180,440)(160,440)
	(160,420)(180,420)
\path(200,420)(200,440)(180,440)
	(180,420)(200,420)
\path(200,400)(200,420)(180,420)
	(180,400)(200,400)
\path(180,400)(180,420)(160,420)
	(160,400)(180,400)
\path(260,340)(260,360)(240,360)
	(240,340)(260,340)
\path(280,340)(280,360)(260,360)
	(260,340)(280,340)
\path(280,320)(280,340)(260,340)
	(260,320)(280,320)
\path(260,320)(260,340)(240,340)
	(240,320)(260,320)
\path(300,420)(300,440)(280,440)
	(280,420)(300,420)
\path(320,420)(320,440)(300,440)
	(300,420)(320,420)
\path(320,400)(320,420)(300,420)
	(300,400)(320,400)
\path(300,400)(300,420)(280,420)
	(280,400)(300,400)
\path(460,260)(460,280)(440,280)
	(440,260)(460,260)
\path(480,260)(480,280)(460,280)
	(460,260)(480,260)
\path(480,240)(480,260)(460,260)
	(460,240)(480,240)
\path(460,240)(460,260)(440,260)
	(440,240)(460,240)
\path(500,340)(500,360)(480,360)
	(480,340)(500,340)
\path(520,340)(520,360)(500,360)
	(500,340)(520,340)
\path(520,320)(520,340)(500,340)
	(500,320)(520,320)
\path(500,320)(500,340)(480,340)
	(480,320)(500,320)
\path(580,260)(580,280)(560,280)
	(560,260)(580,260)
\path(600,260)(600,280)(580,280)
	(580,260)(600,260)
\path(600,240)(600,260)(580,260)
	(580,240)(600,240)
\path(580,240)(580,260)(560,260)
	(560,240)(580,240)
\path(380,380)(380,400)(360,400)
	(360,380)(380,380)
\path(400,380)(400,400)(380,400)
	(380,380)(400,380)
\path(400,360)(400,380)(380,380)
	(380,360)(400,360)
\path(380,360)(380,380)(360,380)
	(360,360)(380,360)
\path(380,300)(380,320)(360,320)
	(360,300)(380,300)
\path(400,300)(400,320)(380,320)
	(380,300)(400,300)
\path(400,280)(400,300)(380,300)
	(380,280)(400,280)
\path(380,280)(380,300)(360,300)
	(360,280)(380,280)
\path(420,20)(420,40)(400,40)
	(400,20)(420,20)
\path(440,20)(440,40)(420,40)
	(420,20)(440,20)
\path(420,0)(420,20)(400,20)
	(400,0)(420,0)
\path(440,0)(440,20)(420,20)
	(420,0)(440,0)
\path(340,100)(340,120)(320,120)
	(320,100)(340,100)
\path(360,100)(360,120)(340,120)
	(340,100)(360,100)
\path(360,80)(360,100)(340,100)
	(340,80)(360,80)
\path(340,80)(340,100)(320,100)
	(320,80)(340,80)
\path(260,180)(260,200)(240,200)
	(240,180)(260,180)
\path(280,180)(280,200)(260,200)
	(260,180)(280,180)
\path(280,160)(280,180)(260,180)
	(260,160)(280,160)
\path(260,160)(260,180)(240,180)
	(240,160)(260,160)
\path(220,100)(220,120)(200,120)
	(200,100)(220,100)
\path(240,100)(240,120)(220,120)
	(220,100)(240,100)
\path(240,80)(240,100)(220,100)
	(220,80)(240,80)
\path(220,80)(220,100)(200,100)
	(200,80)(220,80)
\path(140,180)(140,200)(120,200)
	(120,180)(140,180)
\path(160,180)(160,200)(140,200)
	(140,180)(160,180)
\path(160,160)(160,180)(140,180)
	(140,160)(160,160)
\path(140,160)(140,180)(120,180)
	(120,160)(140,160)
\path(20,180)(20,200)(0,200)
	(0,180)(20,180)
\path(40,180)(40,200)(20,200)
	(20,180)(40,180)
\path(40,160)(40,180)(20,180)
	(20,160)(40,160)
\path(20,160)(20,180)(0,180)
	(0,160)(20,160)
\path(45,180)(115,180)
\path(107.000,178.000)(115.000,180.000)(107.000,182.000)
\path(165,180)(235,180)
\path(227.000,178.000)(235.000,180.000)(227.000,182.000)
\path(245,100)(315,100)
\path(307.000,98.000)(315.000,100.000)(307.000,102.000)
\path(285,155)(315,125)
\path(307.929,129.243)(315.000,125.000)(310.757,132.071)
\path(165,155)(195,125)
\path(187.929,129.243)(195.000,125.000)(190.757,132.071)
\path(365,75)(395,45)
\path(387.929,49.243)(395.000,45.000)(390.757,52.071)
\path(340,660)(340,680)(320,680)
	(320,660)(340,660)
\path(360,660)(360,680)(340,680)
	(340,660)(360,660)
\path(360,640)(360,660)(340,660)
	(340,640)(360,640)
\path(340,640)(340,660)(320,660)
	(320,640)(340,640)
\path(420,580)(420,600)(400,600)
	(400,580)(420,580)
\path(440,580)(440,600)(420,600)
	(420,580)(440,580)
\path(440,560)(440,580)(420,580)
	(420,560)(440,560)
\path(420,560)(420,580)(400,580)
	(400,560)(420,560)
\path(540,580)(540,600)(520,600)
	(520,580)(540,580)
\path(560,580)(560,600)(540,600)
	(540,580)(560,580)
\path(560,560)(560,580)(540,580)
	(540,560)(560,560)
\path(540,560)(540,580)(520,580)
	(520,560)(540,560)
\path(500,500)(500,520)(480,520)
	(480,500)(500,500)
\path(520,500)(520,520)(500,520)
	(500,500)(520,500)
\path(520,480)(520,500)(500,500)
	(500,480)(520,480)
\path(500,480)(500,500)(480,500)
	(480,480)(500,480)
\path(620,500)(620,520)(600,520)
	(600,500)(620,500)
\path(640,500)(640,520)(620,520)
	(620,500)(640,500)
\path(640,480)(640,500)(620,500)
	(620,480)(640,480)
\path(620,480)(620,500)(600,500)
	(600,480)(620,480)
\path(740,500)(740,520)(720,520)
	(720,500)(740,500)
\path(760,500)(760,520)(740,520)
	(740,500)(760,500)
\path(760,480)(760,500)(740,500)
	(740,480)(760,480)
\path(740,480)(740,500)(720,500)
	(720,480)(740,480)
\path(445,580)(515,580)
\path(507.000,578.000)(515.000,580.000)(507.000,582.000)
\path(525,500)(595,500)
\path(587.000,498.000)(595.000,500.000)(587.000,502.000)
\path(645,500)(715,500)
\path(707.000,498.000)(715.000,500.000)(707.000,502.000)
\path(445,555)(475,525)
\path(467.929,529.243)(475.000,525.000)(470.757,532.071)
\path(565,555)(595,525)
\path(587.929,529.243)(595.000,525.000)(590.757,532.071)
\path(365,635)(395,605)
\path(387.929,609.243)(395.000,605.000)(390.757,612.071)
\path(450.105,296.636)(455.000,290.000)(453.803,298.159)
\put(242.500,202.500){\arc{459.619}{5.4978}{5.8926}}
\path(347.704,408.844)(355.000,405.000)(350.372,411.825)
\put(312.500,357.500){\arc{127.475}{4.9098}{5.4423}}
\path(346.816,303.986)(355.000,305.000)(347.302,307.956)
\put(369.167,420.833){\arc{233.393}{1.6925}{2.3764}}
\path(470.689,332.971)(475.000,340.000)(467.888,335.826)
\put(363.750,453.438){\arc{317.772}{0.7951}{1.3082}}
\path(362.352,347.191)(365.000,355.000)(358.988,349.356)
\put(582.500,215.000){\arc{517.325}{2.7862}{3.7135}}
\path(301.155,386.835)(300.000,395.000)(297.177,387.252)
\put(-638.333,493.333){\arc{1886.943}{0.1044}{0.3106}}
\path(258.730,306.852)(260.000,315.000)(255.046,308.408)
\put(728.571,117.143){\arc{1017.263}{3.1570}{3.5411}}
\path(541.155,546.835)(540.000,555.000)(537.177,547.252)
\put(-398.333,653.333){\arc{1886.943}{0.1044}{0.3106}}
\path(497.687,467.085)(500.000,475.000)(494.235,469.104)
\put(765.000,320.000){\arc{614.003}{3.0273}{3.6708}}
\path(205,420)(275,420)
\path(267.000,418.000)(275.000,420.000)(267.000,422.000)
\path(485,260)(555,260)
\path(547.000,258.000)(555.000,260.000)(547.000,262.000)
\path(205,395)(235,365)
\path(227.929,369.243)(235.000,365.000)(230.757,372.071)
\path(525,315)(555,285)
\path(547.929,289.243)(555.000,285.000)(550.757,292.071)
\path(140,205)(180,395)
\path(180.309,386.760)(180.000,395.000)(176.395,387.584)
\path(300,445)(340,635)
\path(340.309,626.760)(340.000,635.000)(336.395,627.584)
\path(580,285)(620,475)
\path(620.309,466.760)(620.000,475.000)(616.395,467.584)
\path(420,45)(460,235)
\path(460.309,226.760)(460.000,235.000)(456.395,227.584)
\path(380,405)(420,555)
\path(419.871,546.755)(420.000,555.000)(416.006,547.785)
\put(75,185){\makebox(0,0)[lb]{\smash{{{\SetFigFont{12}{14.4}{bf}2}}}}}
\put(195,185){\makebox(0,0)[lb]{\smash{{{\SetFigFont{12}{14.4}{bf}2}}}}}
\put(275,105){\makebox(0,0)[lb]{\smash{{{\SetFigFont{12}{14.4}{bf}2}}}}}
\put(235,425){\makebox(0,0)[lb]{\smash{{{\SetFigFont{12}{14.4}{bf}2}}}}}
\put(315,325){\makebox(0,0)[lb]{\smash{{{\SetFigFont{12}{14.4}{bf}2}}}}}
\put(420,310){\makebox(0,0)[lb]{\smash{{{\SetFigFont{12}{14.4}{bf}2}}}}}
\put(515,265){\makebox(0,0)[lb]{\smash{{{\SetFigFont{12}{14.4}{bf}2}}}}}
\put(475,585){\makebox(0,0)[lb]{\smash{{{\SetFigFont{12}{14.4}{bf}2}}}}}
\put(555,505){\makebox(0,0)[lb]{\smash{{{\SetFigFont{12}{14.4}{bf}2}}}}}
\put(675,505){\makebox(0,0)[lb]{\smash{{{\SetFigFont{12}{14.4}{bf}2}}}}}
\put(180,145){\makebox(0,0)[lb]{\smash{{{\SetFigFont{12}{14.4}{bf}1}}}}}
\put(300,145){\makebox(0,0)[lb]{\smash{{{\SetFigFont{12}{14.4}{bf}1}}}}}
\put(380,65){\makebox(0,0)[lb]{\smash{{{\SetFigFont{12}{14.4}{bf}1}}}}}
\put(380,625){\makebox(0,0)[lb]{\smash{{{\SetFigFont{12}{14.4}{bf}1}}}}}
\put(460,545){\makebox(0,0)[lb]{\smash{{{\SetFigFont{12}{14.4}{bf}1}}}}}
\put(580,545){\makebox(0,0)[lb]{\smash{{{\SetFigFont{12}{14.4}{bf}1}}}}}
\put(220,385){\makebox(0,0)[lb]{\smash{{{\SetFigFont{12}{14.4}{bf}1}}}}}
\put(540,305){\makebox(0,0)[lb]{\smash{{{\SetFigFont{12}{14.4}{bf}1}}}}}
\put(340,420){\makebox(0,0)[lb]{\smash{{{\SetFigFont{12}{14.4}{bf}1}}}}}
\put(425,350){\makebox(0,0)[lb]{\smash{{{\SetFigFont{12}{14.4}{bf}1}}}}}
\put(145,300){\makebox(0,0)[lb]{\smash{{{\SetFigFont{12}{14.4}{bf}3}}}}}
\put(220,225){\makebox(0,0)[lb]{\smash{{{\SetFigFont{12}{14.4}{bf}3}}}}}
\put(285,270){\makebox(0,0)[lb]{\smash{{{\SetFigFont{12}{14.4}{bf}3}}}}}
\put(310,225){\makebox(0,0)[lb]{\smash{{{\SetFigFont{12}{14.4}{bf}3}}}}}
\put(430,145){\makebox(0,0)[lb]{\smash{{{\SetFigFont{12}{14.4}{bf}3}}}}}
\put(590,380){\makebox(0,0)[lb]{\smash{{{\SetFigFont{12}{14.4}{bf}3}}}}}
\put(525,445){\makebox(0,0)[lb]{\smash{{{\SetFigFont{12}{14.4}{bf}3}}}}}
\put(455,390){\makebox(0,0)[lb]{\smash{{{\SetFigFont{12}{14.4}{bf}3}}}}}
\put(305,535){\makebox(0,0)[lb]{\smash{{{\SetFigFont{12}{14.4}{bf}3}}}}}
\put(385,475){\makebox(0,0)[lb]{\smash{{{\SetFigFont{12}{14.4}{bf}3}}}}}
\put(5,185){\makebox(0,0)[lb]{\smash{{{\SetFigFont{12}{14.4}{bf}2}}}}}
\put(25,185){\makebox(0,0)[lb]{\smash{{{\SetFigFont{12}{14.4}{bf}2}}}}}
\put(5,165){\makebox(0,0)[lb]{\smash{{{\SetFigFont{12}{14.4}{bf}1}}}}}
\put(25,165){\makebox(0,0)[lb]{\smash{{{\SetFigFont{12}{14.4}{bf}1}}}}}
\put(125,185){\makebox(0,0)[lb]{\smash{{{\SetFigFont{12}{14.4}{bf}2}}}}}
\put(145,185){\makebox(0,0)[lb]{\smash{{{\SetFigFont{12}{14.4}{bf}3}}}}}
\put(125,165){\makebox(0,0)[lb]{\smash{{{\SetFigFont{12}{14.4}{bf}1}}}}}
\put(145,165){\makebox(0,0)[lb]{\smash{{{\SetFigFont{12}{14.4}{bf}1}}}}}
\put(205,105){\makebox(0,0)[lb]{\smash{{{\SetFigFont{12}{14.4}{bf}2}}}}}
\put(225,105){\makebox(0,0)[lb]{\smash{{{\SetFigFont{12}{14.4}{bf}3}}}}}
\put(205,85){\makebox(0,0)[lb]{\smash{{{\SetFigFont{12}{14.4}{bf}1}}}}}
\put(225,85){\makebox(0,0)[lb]{\smash{{{\SetFigFont{12}{14.4}{bf}2}}}}}
\put(325,105){\makebox(0,0)[lb]{\smash{{{\SetFigFont{12}{14.4}{bf}3}}}}}
\put(345,105){\makebox(0,0)[lb]{\smash{{{\SetFigFont{12}{14.4}{bf}3}}}}}
\put(325,85){\makebox(0,0)[lb]{\smash{{{\SetFigFont{12}{14.4}{bf}1}}}}}
\put(345,85){\makebox(0,0)[lb]{\smash{{{\SetFigFont{12}{14.4}{bf}2}}}}}
\put(405,25){\makebox(0,0)[lb]{\smash{{{\SetFigFont{12}{14.4}{bf}3}}}}}
\put(425,25){\makebox(0,0)[lb]{\smash{{{\SetFigFont{12}{14.4}{bf}3}}}}}
\put(405,5){\makebox(0,0)[lb]{\smash{{{\SetFigFont{12}{14.4}{bf}2}}}}}
\put(425,5){\makebox(0,0)[lb]{\smash{{{\SetFigFont{12}{14.4}{bf}2}}}}}
\put(245,185){\makebox(0,0)[lb]{\smash{{{\SetFigFont{12}{14.4}{bf}3}}}}}
\put(265,185){\makebox(0,0)[lb]{\smash{{{\SetFigFont{12}{14.4}{bf}3}}}}}
\put(245,165){\makebox(0,0)[lb]{\smash{{{\SetFigFont{12}{14.4}{bf}1}}}}}
\put(265,165){\makebox(0,0)[lb]{\smash{{{\SetFigFont{12}{14.4}{bf}1}}}}}
\put(165,425){\makebox(0,0)[lb]{\smash{{{\SetFigFont{12}{14.4}{bf}2}}}}}
\put(185,425){\makebox(0,0)[lb]{\smash{{{\SetFigFont{12}{14.4}{bf}4}}}}}
\put(165,405){\makebox(0,0)[lb]{\smash{{{\SetFigFont{12}{14.4}{bf}1}}}}}
\put(185,405){\makebox(0,0)[lb]{\smash{{{\SetFigFont{12}{14.4}{bf}1}}}}}
\put(285,425){\makebox(0,0)[lb]{\smash{{{\SetFigFont{12}{14.4}{bf}3}}}}}
\put(305,425){\makebox(0,0)[lb]{\smash{{{\SetFigFont{12}{14.4}{bf}4}}}}}
\put(285,405){\makebox(0,0)[lb]{\smash{{{\SetFigFont{12}{14.4}{bf}1}}}}}
\put(305,405){\makebox(0,0)[lb]{\smash{{{\SetFigFont{12}{14.4}{bf}1}}}}}
\put(245,345){\makebox(0,0)[lb]{\smash{{{\SetFigFont{12}{14.4}{bf}2}}}}}
\put(265,345){\makebox(0,0)[lb]{\smash{{{\SetFigFont{12}{14.4}{bf}4}}}}}
\put(245,325){\makebox(0,0)[lb]{\smash{{{\SetFigFont{12}{14.4}{bf}1}}}}}
\put(265,325){\makebox(0,0)[lb]{\smash{{{\SetFigFont{12}{14.4}{bf}2}}}}}
\put(365,385){\makebox(0,0)[lb]{\smash{{{\SetFigFont{12}{14.4}{bf}3}}}}}
\put(385,385){\makebox(0,0)[lb]{\smash{{{\SetFigFont{12}{14.4}{bf}4}}}}}
\put(365,365){\makebox(0,0)[lb]{\smash{{{\SetFigFont{12}{14.4}{bf}1}}}}}
\put(385,365){\makebox(0,0)[lb]{\smash{{{\SetFigFont{12}{14.4}{bf}2}}}}}
\put(365,305){\makebox(0,0)[lb]{\smash{{{\SetFigFont{12}{14.4}{bf}2}}}}}
\put(385,305){\makebox(0,0)[lb]{\smash{{{\SetFigFont{12}{14.4}{bf}4}}}}}
\put(365,285){\makebox(0,0)[lb]{\smash{{{\SetFigFont{12}{14.4}{bf}1}}}}}
\put(385,285){\makebox(0,0)[lb]{\smash{{{\SetFigFont{12}{14.4}{bf}3}}}}}
\put(445,265){\makebox(0,0)[lb]{\smash{{{\SetFigFont{12}{14.4}{bf}3}}}}}
\put(465,265){\makebox(0,0)[lb]{\smash{{{\SetFigFont{12}{14.4}{bf}4}}}}}
\put(445,245){\makebox(0,0)[lb]{\smash{{{\SetFigFont{12}{14.4}{bf}2}}}}}
\put(465,245){\makebox(0,0)[lb]{\smash{{{\SetFigFont{12}{14.4}{bf}2}}}}}
\put(485,345){\makebox(0,0)[lb]{\smash{{{\SetFigFont{12}{14.4}{bf}3}}}}}
\put(505,345){\makebox(0,0)[lb]{\smash{{{\SetFigFont{12}{14.4}{bf}4}}}}}
\put(485,325){\makebox(0,0)[lb]{\smash{{{\SetFigFont{12}{14.4}{bf}1}}}}}
\put(505,325){\makebox(0,0)[lb]{\smash{{{\SetFigFont{12}{14.4}{bf}3}}}}}
\put(565,265){\makebox(0,0)[lb]{\smash{{{\SetFigFont{12}{14.4}{bf}3}}}}}
\put(585,265){\makebox(0,0)[lb]{\smash{{{\SetFigFont{12}{14.4}{bf}4}}}}}
\put(565,245){\makebox(0,0)[lb]{\smash{{{\SetFigFont{12}{14.4}{bf}2}}}}}
\put(585,245){\makebox(0,0)[lb]{\smash{{{\SetFigFont{12}{14.4}{bf}3}}}}}
\put(325,665){\makebox(0,0)[lb]{\smash{{{\SetFigFont{12}{14.4}{bf}4}}}}}
\put(345,665){\makebox(0,0)[lb]{\smash{{{\SetFigFont{12}{14.4}{bf}4}}}}}
\put(325,645){\makebox(0,0)[lb]{\smash{{{\SetFigFont{12}{14.4}{bf}1}}}}}
\put(345,645){\makebox(0,0)[lb]{\smash{{{\SetFigFont{12}{14.4}{bf}1}}}}}
\put(405,585){\makebox(0,0)[lb]{\smash{{{\SetFigFont{12}{14.4}{bf}4}}}}}
\put(425,585){\makebox(0,0)[lb]{\smash{{{\SetFigFont{12}{14.4}{bf}4}}}}}
\put(405,565){\makebox(0,0)[lb]{\smash{{{\SetFigFont{12}{14.4}{bf}1}}}}}
\put(425,565){\makebox(0,0)[lb]{\smash{{{\SetFigFont{12}{14.4}{bf}2}}}}}
\put(525,585){\makebox(0,0)[lb]{\smash{{{\SetFigFont{12}{14.4}{bf}4}}}}}
\put(545,585){\makebox(0,0)[lb]{\smash{{{\SetFigFont{12}{14.4}{bf}4}}}}}
\put(525,565){\makebox(0,0)[lb]{\smash{{{\SetFigFont{12}{14.4}{bf}1}}}}}
\put(545,565){\makebox(0,0)[lb]{\smash{{{\SetFigFont{12}{14.4}{bf}3}}}}}
\put(485,505){\makebox(0,0)[lb]{\smash{{{\SetFigFont{12}{14.4}{bf}4}}}}}
\put(505,505){\makebox(0,0)[lb]{\smash{{{\SetFigFont{12}{14.4}{bf}4}}}}}
\put(505,485){\makebox(0,0)[lb]{\smash{{{\SetFigFont{12}{14.4}{bf}2}}}}}
\put(485,485){\makebox(0,0)[lb]{\smash{{{\SetFigFont{12}{14.4}{bf}2}}}}}
\put(605,505){\makebox(0,0)[lb]{\smash{{{\SetFigFont{12}{14.4}{bf}4}}}}}
\put(625,505){\makebox(0,0)[lb]{\smash{{{\SetFigFont{12}{14.4}{bf}4}}}}}
\put(605,485){\makebox(0,0)[lb]{\smash{{{\SetFigFont{12}{14.4}{bf}2}}}}}
\put(625,485){\makebox(0,0)[lb]{\smash{{{\SetFigFont{12}{14.4}{bf}3}}}}}
\put(725,505){\makebox(0,0)[lb]{\smash{{{\SetFigFont{12}{14.4}{bf}4}}}}}
\put(745,505){\makebox(0,0)[lb]{\smash{{{\SetFigFont{12}{14.4}{bf}4}}}}}
\put(725,485){\makebox(0,0)[lb]{\smash{{{\SetFigFont{12}{14.4}{bf}3}}}}}
\put(745,485){\makebox(0,0)[lb]{\smash{{{\SetFigFont{12}{14.4}{bf}3}}}}}
\end{picture}
\caption{Crystal graph of the $U_q(\sl_4)$-module $V_{(2,\,2)}$\label{G22}}
\end{center}
\end{figure}

If one restricts the crystal graph $\Gamma_{\lambda}$ to its edges of colour
$i$,
one obtains a decomposition of this graph into {\it strings of colour $i$}.
Each tableau $t$ belongs to a unique string of colour $i$, possibly reduced
to $t$.

The graph $\Gamma_{\lambda}$ can also be decomposed in a different way.
Indeed, the operators $\sigma_i$ satisfy Moore-Coxeter relations \cite{LS1},
which allows to define an action of the Weyl group $W=\S_{n+1}$ on the crystal
graph $\Gamma_{\lambda}$ and to decompose it into {\it orbits}.
We shall denote by $\O_t$ the orbit of the tableau $t$ under $\S_{n+1}$.

Note that a similar action of the Weyl group on the crystal graph can be
defined for the other root systems. In geometric terms, the
$i$-th generator of the Weyl group sends a vertex to its
mirror image through the middle of its string of colour $i$.
For example, on the graph of Figure~\ref{G22} one can see four orbits.
Reading tableaux columnwise, these are the two fixed points $2143$
and $3142$, a $6$-element orbit $\{2121,\,3131,\,4141,\,3232,\,4242,\,4343\}$,
and a $12$-element orbit constituted of the remaining tableaux.

As in the classical case, the multiplicity of the
weight $\mu$ in $V_\lambda$ is equal to the number $K_{\lambda\,\mu}$
of tableaux of shape $\lambda$ and weight $\mu$, that is, to the number
of orbits $\O_t$ parametrized by these tableaux.

By definition, the charge is constant along the orbits $\O_t$,
and can therefore be regarded as a statistic on the set of these orbits.

We shall now introduce a different statistic $d(t)$ on tableaux, which
reflects the geometry of the crystal graph around $t$. Denote by
$d_i(t)$ the {\it exponent of $t$ in direction $i$}, defined by
$$
d_i(t) = {\rm min}\{{\rm max}\{k \mid \epsilon_i^k(t) \not = 0\},\,
					{\rm max}\{k \mid \phi_i^k(t) \not = 0\}\} \,.
$$
In other words, $d_i(t)$ is the distance from $t$ to the nearest end
of its string of colour $i$. Define then
$$
d(t) = \sum _{i=1}^n i\,d_i(t) \,.
$$
The statistic $d(t')$ is not constant for $t'$ in the orbit $\O_t$. But
one has
\begin{theorem}
{\rm (i)} For any tableau $t$, the arithmetic mean of the integers
$d(t'),\,t'\in \O_t$, is an integer denoted by $b(t)$.

{\rm (ii)} $b(t)$ is equal to the charge $c(u)$ of the image $u$ of $t$
under the involution $\Omega_2$.

{\rm (iii)} The Kostka-Foulkes polynomial $K_{\lambda\,\mu}(q)$ is equal to
$$
K_{\lambda\,\mu}(q) = \sum_{t\in \tab(\lambda,\,\mu)} q^{b(t)} \,.
$$
\end{theorem}

\begin{example}{\rm
We take $n=3$, $\lambda = (3,\,2)$, $\mu = (2,\,1,\,1,\,1)$. There are three
tableaux
in $\tab(\lambda,\,\mu)$

\setlength{\unitlength}{0.25pt}
\centerline{
\begin{picture}(1000,130)
\put(0,0){\begin{picture}(250,125)
\put(100,0){\begin{picture}(150,100)
\put(0,0){\framebox(50,50){1}}
\put(50,0){\framebox(50,50){1}}
\put(100,0){\framebox(50,50){2}}
\put(0,50){\framebox(50,50){3}}
\put(50,50){\framebox(50,50){4}}
\end{picture}}
\put(0,20){$t =$}
\end{picture}}
\put(325,0){\begin{picture}(250,125)
\put(100,0){\begin{picture}(150,100)
\put(0,0){\framebox(50,50){1}}
\put(50,0){\framebox(50,50){1}}
\put(100,0){\framebox(50,50){3}}
\put(0,50){\framebox(50,50){2}}
\put(50,50){\framebox(50,50){4}}
\end{picture}}
\put(0,20){$u =$}
\end{picture}}
\put(650,0){\begin{picture}(250,125)
\put(100,0){\begin{picture}(150,100)
\put(0,0){\framebox(50,50){1}}
\put(50,0){\framebox(50,50){1}}
\put(100,0){\framebox(50,50){4}}
\put(0,50){\framebox(50,50){2}}
\put(50,50){\framebox(50,50){3}}
\end{picture}}
\put(0,20){$v =$}
\end{picture}}
\end{picture}}

\medskip
\noindent
The orbit of $t$ is made of the following tableaux $t'$

\centerline{
\begin{picture}(900,150)
\put(0,0){\begin{picture}(150,100)
\put(0,0){\framebox(50,50){1}}
\put(50,0){\framebox(50,50){1}}
\put(100,0){\framebox(50,50){2}}
\put(0,50){\framebox(50,50){3}}
\put(50,50){\framebox(50,50){4}}
                  \end{picture}}
\put(200,0){\begin{picture}(150,100)
\put(0,0){\framebox(50,50){1}}
\put(50,0){\framebox(50,50){2}}
\put(100,0){\framebox(50,50){2}}
\put(0,50){\framebox(50,50){3}}
\put(50,50){\framebox(50,50){4}}
                  \end{picture}}
\put(400,0){\begin{picture}(150,100)
\put(0,0){\framebox(50,50){1}}
\put(50,0){\framebox(50,50){2}}
\put(100,0){\framebox(50,50){3}}
\put(0,50){\framebox(50,50){3}}
\put(50,50){\framebox(50,50){4}}
                  \end{picture}}
\put(600,0){\begin{picture}(150,100)
\put(0,0){\framebox(50,50){1}}
\put(50,0){\framebox(50,50){2}}
\put(100,0){\framebox(50,50){3}}
\put(0,50){\framebox(50,50){4}}
\put(50,50){\framebox(50,50){4}}
                  \end{picture}}
\end{picture}}

\medskip\noindent
whose numbers $d(t')$ are respectively equal to $4,\,4,\,1,\,3$.
Similarly, the orbit $\O_u$
is made of the following tableaux $u'$

\centerline{
\begin{picture}(900,150)
\put(0,0){\begin{picture}(150,100)
\put(0,0){\framebox(50,50){1}}
\put(50,0){\framebox(50,50){1}}
\put(100,0){\framebox(50,50){3}}
\put(0,50){\framebox(50,50){2}}
\put(50,50){\framebox(50,50){4}}
\end{picture}}
\put(200,0){\begin{picture}(150,100)
\put(0,0){\framebox(50,50){1}}
\put(50,0){\framebox(50,50){2}}
\put(100,0){\framebox(50,50){3}}
\put(0,50){\framebox(50,50){2}}
\put(50,50){\framebox(50,50){4}}
\end{picture}}
\put(400,0){\begin{picture}(150,100)
\put(0,0){\framebox(50,50){1}}
\put(50,0){\framebox(50,50){3}}
\put(100,0){\framebox(50,50){3}}
\put(0,50){\framebox(50,50){2}}
\put(50,50){\framebox(50,50){4}}
\end{picture}}
\put(600,0){\begin{picture}(150,100)
\put(0,0){\framebox(50,50){1}}
\put(50,0){\framebox(50,50){3}}
\put(100,0){\framebox(50,50){4}}
\put(0,50){\framebox(50,50){2}}
\put(50,50){\framebox(50,50){4}}
\end{picture}}
\end{picture}}

\medskip\noindent
whose numbers $d(u')$ are all equal to $2$.
Finally, the orbit $\O_v$
is made of the tableaux $v'$

\centerline{
\begin{picture}(900,150)
\put(0,0){\begin{picture}(150,100)
\put(0,0){\framebox(50,50){1}}
\put(50,0){\framebox(50,50){1}}
\put(100,0){\framebox(50,50){4}}
\put(0,50){\framebox(50,50){2}}
\put(50,50){\framebox(50,50){3}}
\end{picture}}
\put(200,0){\begin{picture}(150,100)
\put(0,0){\framebox(50,50){1}}
\put(50,0){\framebox(50,50){2}}
\put(100,0){\framebox(50,50){4}}
\put(0,50){\framebox(50,50){2}}
\put(50,50){\framebox(50,50){3}}
\end{picture}}
\put(400,0){\begin{picture}(150,100)
\put(0,0){\framebox(50,50){1}}
\put(50,0){\framebox(50,50){2}}
\put(100,0){\framebox(50,50){4}}
\put(0,50){\framebox(50,50){3}}
\put(50,50){\framebox(50,50){3}}
\end{picture}}
\put(600,0){\begin{picture}(150,100)
\put(0,0){\framebox(50,50){1}}
\put(50,0){\framebox(50,50){2}}
\put(100,0){\framebox(50,50){4}}
\put(0,50){\framebox(50,50){3}}
\put(50,50){\framebox(50,50){4}}
\end{picture}}
\end{picture}}

\medskip\noindent
whose numbers $d(v')$ are respectively $5,\,3,\,4,\,4$.
Thus, $b(t)=3=c(v)$, $b(u)=2=c(u)$, $b(v)=4=c(t)$, and one
has $\Omega_2(t) \in \O_v$, $\Omega_2(u) \in O_u$, $\Omega_2(v) \in \O_t$.
As a result, we get $K_{\lambda\,\mu}(q) = q^2+q^3+q^4$.
}
\end{example}

\Proof We establish (ii) which clearly implies (i) and (iii). Let
$d'(t) = d(\Omega_2(t)) = \sum_{1\le i \le n} (n-i+1) d_i(t)$ by the
symmetry property of $\Gamma_\lambda$ mentionned in section~\ref{PLACTIC}.
Denote by $b'(t)$ the arithmetic mean of the integers $d'(t')$ for
$t' \in {\cal O}_t$. More generally, given a subset $S$ of ${\cal O}_t$,
we denote by $B'(S)$ the arithmetic mean of the $d'(t')$ for $t' \in S$.

The proof proceeds by induction. One first checks by direct computation that
for
the row tableau $t_\mu$ of weight the partition $\mu$, there holds
\begin{equation}
b'(t_\mu) = \|\mu\| \,.
\end{equation}
Therefore, from the definition of $c(t)$, it is sufficient to prove that
\begin{equation}
{\cal C}(t) = s \quad \Longrightarrow \quad b'(s) = b'(t) + 1 \,,
\end{equation}
where ${\cal C}(t)$ is the initial cyclage of $t$ (see section~\ref{PLACTIC}).
To this end, we decompose the orbit ${\cal O}_t$ into chains in the
following way. For two tableaux $u,\,v \in {\cal O}_t$, write
$u \leadsto v$ if there exists an $i\in \{2,\,\ldots ,\,n+1\}$ such that
$\sigma_i(u)=v$ and $a_{i+1}$ (resp. $a_i$) is the first letter of the
row reading of $u$ (resp. of $v$). The connected components of the resulting
graph are linear graphs called {\it chains}. Now, if $\delta = {\cal
C}(\gamma)$
is the set of tableaux obtained by applying the initial cyclage operator
to a chain $\gamma$ of ${\cal O}_t$, then
\begin{equation}\label{EQUA}
B'(\delta) = B'(\gamma) + 1\,.
\end{equation}
Indeed, it is possible to describe explicitely the difference between the
$n$-uple of exponents $\d(t)=(d_1(t),\ldots ,d_n(t))$ of a tableau $t$ of
the chain, and that of its image ${\cal C}(t)$.
It turns out that this difference vector is zero, except at the ends
of the chain.

More precisely, let $\gamma=\{u_i\leadsto u_{i-1}\leadsto \ldots \leadsto
u_h\}$
be a chain,
where the indices have been chosen such that $a_j$ is the first letter
of $u_j$, and let ${\cal C}(u_j)=v_j$. Note that since $u_i$ is not
a row tableau, one has $h\ge 2$. Let $(\e_j)$ be the canonical basis
of $\Z^{n}$. We also set for convenience $\e_{n+1} = 0$. Then,
the difference vectors are as follows if $i>h$
\begin{eqnarray}
\d(v_i) - \d(u_i) & = & - \e_i \\
\d(v_k) - \d(u_k) & = &  0, \ \ \ i>k>h \\
\d(v_h) - \d(u_h) & = &  \e_{h-1}
\end{eqnarray}
and if $i=h$
\begin{equation}
\d(v_i) - \d(u_i)\  =\  \e_{i-1} - \e_i \,. \ \ \ \
\end{equation}
This clearly implies (\ref{EQUA}). \cqfd

\noindent
{\bf Remark}

\medskip\noindent
As observed by Terada \cite{Ter}, the standard polynomials
$K_{\lambda,(1^n)}(q)$
can also be interpreted in terms of the Kazhdan-Lusztig basis of the
irreducible $S_n$-module $W_\lambda$.

\section{Refinement of the generalized exponents}

The polynomials $K_{\lambda,(k^{n+1})}(q)$ are of particular importance
for representation theory. The coefficient of $q^d$ in
$K_{\lambda,(k^{n+1})}(q)$
is equal to the multiplicity of the irreducible $SL_{n+1}$-module
$V_\lambda$ in the homogeneous component of degree $d$ of the
affine coordinate ring of the variety of nilpotent matrices in
$\sl_{n+1}$ \cite{He}, or equivalently of the space of $SL_{n+1}$-harmonic
polynomials \cite{Gupta1}. Similar polynomials can be defined for
any reductive complex algebraic groups \cite{Ko}. The exponents of
the nonzero terms of  $K_{\lambda,(k^{n+1})}(q)$ have been called
by Kostant the {\it generalized exponents} of the module $V_\lambda$.

In the case where $\mu$ is a partition of rectangular shape $\mu=(k^{n+1})$,
the previous construction becomes simpler
because all the elements of $\tab(\lambda,(k^{n+1}))$ remain fixed
under $\S_{n+1}$. One can then define polynomials in
several variables
$$
\K_{\lambda,(k^{n+1})}(x_1,\ldots,x_n) =
\sum_{t\in \tab(\lambda,(k^{n+1}))}\prod_{i=1}^n x_i^{d_i(t)}
$$
such that $K_{\lambda,(k^{n+1})}(q)=\K_{\lambda,(k^{n+1})}(q,q^2,\ldots,q^n)$.

\begin{example}{\rm
{\rm (i)} $\K_{(3,3,2),(2,2,2,2)}(x_1,\,x_2,\,x_3) =
x_1x_2 + x_1x_3 + x_2x_3$.

{\rm (ii)} $\K_{(5,2,1),(2,2,2,2)}(x_1,\,x_2,\,x_3) =
x_2x_3^2 + x_2^2x_3 +x_1x_3^2 + 2\,x_1x_2x_3+x_1x_2^2+x_1^2x_3+x_1^2x_2$.

{\rm (iii)} $\K_{(4,3,1),(2,2,2,2)}(x_1,\,x_2,\,x_3) =
x_1x_2+x_2x_3+x_2^2+x_1x_3^2+x_1^2x_3+2\,x_1x_2x_3$.}
\end{example}

These polynomials admit a rather simple description in terms
of plactic operations.

\begin{theorem}
Let $t\in\tab(\lambda,(k^{n+1}))$. Then,

{\rm (i)} There exists a unique tableau $u$ of minimal weight
such that in the plactic monoid, $t\,.\,u$ is a Yamanouchi tableau.
Moreover, $u$ is itself a Yamanouchi tableau and if $y$ is the Yamanouchi
tableau of weight $(k^{n+1})$, one has
$$ t\,u \equiv u\,y \equiv y\,u \ . $$

{\rm (ii)} Denote by $\nu(t)$ the weight of the tableau $u$ defined
in {\rm (i)}, and for a partition $\mu$, set $x_\mu=
x_{\mu_1}x_{\mu_2}\cdots x_{\mu_r}$. Then,
$$
\K_{\lambda,(k^{n+1})}(x_1,\ldots,x_n) =
\sum_{t\in\tab(\lambda,(k^{n+1}))} x_{\nu'(t)} \ ,
$$
where $\nu'(t)$ denotes the conjugate of the partition $\nu(t)$.
\end{theorem}

\begin{example}{\rm With $\lambda=(4,2)$ and $\mu=(2,2,2)$, one has the
following equalities in the plactic monoid

\begin{center}
\setlength{\unitlength}{0.0085in}
\begingroup\makeatletter\ifx\SetFigFont\undefined
\def\x#1#2#3#4#5#6#7\relax{\def\x{#1#2#3#4#5#6}}%
\expandafter\x\fmtname xxxxxx\relax \def\y{splain}%
\ifx\x\y   
\gdef\SetFigFont#1#2#3{%
  \ifnum #1<17\tiny\else \ifnum #1<20\small\else
  \ifnum #1<24\normalsize\else \ifnum #1<29\large\else
  \ifnum #1<34\Large\else \ifnum #1<41\LARGE\else
     \huge\fi\fi\fi\fi\fi\fi
  \csname #3\endcsname}%
\else
\gdef\SetFigFont#1#2#3{\begingroup
  \count@#1\relax \ifnum 25<\count@\count@25\fi
  \def\x{\endgroup\@setsize\SetFigFont{#2pt}}%
  \expandafter\x
    \csname \romannumeral\the\count@ pt\expandafter\endcsname
    \csname @\romannumeral\the\count@ pt\endcsname
  \csname #3\endcsname}%
\fi
\fi\endgroup
\begin{picture}(569,197)(0,-10)
\path(195,160)(195,180)(175,180)
	(175,160)(195,160)
\path(195,140)(195,160)(175,160)
	(175,140)(195,140)
\path(215,140)(215,160)(195,160)
	(195,140)(215,140)
\path(195,120)(195,140)(175,140)
	(175,120)(195,120)
\path(215,120)(215,140)(195,140)
	(195,120)(215,120)
\path(235,120)(235,140)(215,140)
	(215,120)(235,120)
\path(215,160)(215,180)(195,180)
	(195,160)(215,160)
\path(255,120)(255,140)(235,140)
	(235,120)(255,120)
\put(180,165){\makebox(0,0)[lb]{\smash{{{\SetFigFont{12}{14.4}{rm}3}}}}}
\put(200,165){\makebox(0,0)[lb]{\smash{{{\SetFigFont{12}{14.4}{rm}3}}}}}
\put(180,145){\makebox(0,0)[lb]{\smash{{{\SetFigFont{12}{14.4}{rm}2}}}}}
\put(200,145){\makebox(0,0)[lb]{\smash{{{\SetFigFont{12}{14.4}{rm}2}}}}}
\put(180,125){\makebox(0,0)[lb]{\smash{{{\SetFigFont{12}{14.4}{rm}1}}}}}
\put(200,125){\makebox(0,0)[lb]{\smash{{{\SetFigFont{12}{14.4}{rm}1}}}}}
\put(220,125){\makebox(0,0)[lb]{\smash{{{\SetFigFont{12}{14.4}{rm}1}}}}}
\put(240,125){\makebox(0,0)[lb]{\smash{{{\SetFigFont{12}{14.4}{rm}1}}}}}
\put(265,120){\makebox(0,0)[lb]{\smash{{{\SetFigFont{12}{14.4}{rm},}}}}}
\path(20,140)(20,160)(0,160)
	(0,140)(20,140)
\path(20,120)(20,140)(0,140)
	(0,120)(20,120)
\path(40,120)(40,140)(20,140)
	(20,120)(40,120)
\path(40,140)(40,140)(40,140)
	(40,140)(40,140)
\path(60,120)(60,140)(40,140)
	(40,120)(60,120)
\path(40,140)(40,160)(20,160)
	(20,140)(40,140)
\path(80,120)(80,140)(60,140)
	(60,120)(80,120)
\path(120,120)(120,140)(100,140)
	(100,120)(120,120)
\path(140,120)(140,140)(120,140)
	(120,120)(140,120)
\put(5,145){\makebox(0,0)[lb]{\smash{{{\SetFigFont{12}{14.4}{rm}3}}}}}
\put(25,145){\makebox(0,0)[lb]{\smash{{{\SetFigFont{12}{14.4}{rm}3}}}}}
\put(5,125){\makebox(0,0)[lb]{\smash{{{\SetFigFont{12}{14.4}{rm}1}}}}}
\put(25,125){\makebox(0,0)[lb]{\smash{{{\SetFigFont{12}{14.4}{rm}1}}}}}
\put(45,125){\makebox(0,0)[lb]{\smash{{{\SetFigFont{12}{14.4}{rm}2}}}}}
\put(65,125){\makebox(0,0)[lb]{\smash{{{\SetFigFont{12}{14.4}{rm}2}}}}}
\put(105,125){\makebox(0,0)[lb]{\smash{{{\SetFigFont{12}{14.4}{rm}1}}}}}
\put(125,125){\makebox(0,0)[lb]{\smash{{{\SetFigFont{12}{14.4}{rm}1}}}}}
\put(90,125){\makebox(0,0)[lb]{\smash{{{\SetFigFont{12}{14.4}{rm}.}}}}}
\put(155,125){\makebox(0,0)[lb]{\smash{{{\SetFigFont{12}{14.4}{rm}=}}}}}
\path(495,160)(495,180)(475,180)
	(475,160)(495,160)
\path(515,160)(515,180)(495,180)
	(495,160)(515,160)
\path(495,140)(495,160)(475,160)
	(475,140)(495,140)
\path(515,140)(515,160)(495,160)
	(495,140)(515,140)
\path(535,140)(535,160)(515,160)
	(515,140)(535,140)
\path(495,120)(495,140)(475,140)
	(475,120)(495,120)
\path(515,120)(515,140)(495,140)
	(495,120)(515,120)
\path(515,120)(515,140)(515,140)
	(515,120)(515,120)
\path(535,120)(535,140)(515,140)
	(515,120)(535,120)
\path(555,120)(555,140)(535,140)
	(535,120)(555,120)
\put(480,165){\makebox(0,0)[lb]{\smash{{{\SetFigFont{12}{14.4}{rm}3}}}}}
\put(500,165){\makebox(0,0)[lb]{\smash{{{\SetFigFont{12}{14.4}{rm}3}}}}}
\put(480,145){\makebox(0,0)[lb]{\smash{{{\SetFigFont{12}{14.4}{rm}2}}}}}
\put(500,145){\makebox(0,0)[lb]{\smash{{{\SetFigFont{12}{14.4}{rm}2}}}}}
\put(520,145){\makebox(0,0)[lb]{\smash{{{\SetFigFont{12}{14.4}{rm}2}}}}}
\put(480,125){\makebox(0,0)[lb]{\smash{{{\SetFigFont{12}{14.4}{rm}1}}}}}
\put(500,125){\makebox(0,0)[lb]{\smash{{{\SetFigFont{12}{14.4}{rm}1}}}}}
\put(520,125){\makebox(0,0)[lb]{\smash{{{\SetFigFont{12}{14.4}{rm}1}}}}}
\put(540,125){\makebox(0,0)[lb]{\smash{{{\SetFigFont{12}{14.4}{rm}1}}}}}
\put(565,120){\makebox(0,0)[lb]{\smash{{{\SetFigFont{12}{14.4}{rm},}}}}}
\path(320,140)(320,160)(300,160)
	(300,140)(320,140)
\path(320,120)(320,140)(300,140)
	(300,120)(320,120)
\path(340,140)(340,160)(320,160)
	(320,140)(340,140)
\path(340,140)(340,140)(320,140)
	(320,140)(340,140)
\path(340,120)(340,140)(320,140)
	(320,120)(340,120)
\path(360,120)(360,140)(340,140)
	(340,120)(360,120)
\path(380,120)(380,140)(360,140)
	(360,120)(380,120)
\path(420,140)(420,160)(400,160)
	(400,140)(420,140)
\path(440,120)(440,140)(420,140)
	(420,120)(440,120)
\path(420,120)(420,140)(400,140)
	(400,120)(420,120)
\put(305,145){\makebox(0,0)[lb]{\smash{{{\SetFigFont{12}{14.4}{rm}2}}}}}
\put(325,145){\makebox(0,0)[lb]{\smash{{{\SetFigFont{12}{14.4}{rm}3}}}}}
\put(305,125){\makebox(0,0)[lb]{\smash{{{\SetFigFont{12}{14.4}{rm}1}}}}}
\put(325,125){\makebox(0,0)[lb]{\smash{{{\SetFigFont{12}{14.4}{rm}1}}}}}
\put(345,125){\makebox(0,0)[lb]{\smash{{{\SetFigFont{12}{14.4}{rm}2}}}}}
\put(365,125){\makebox(0,0)[lb]{\smash{{{\SetFigFont{12}{14.4}{rm}3}}}}}
\put(405,145){\makebox(0,0)[lb]{\smash{{{\SetFigFont{12}{14.4}{rm}2}}}}}
\put(405,125){\makebox(0,0)[lb]{\smash{{{\SetFigFont{12}{14.4}{rm}1}}}}}
\put(425,125){\makebox(0,0)[lb]{\smash{{{\SetFigFont{12}{14.4}{rm}1}}}}}
\put(390,125){\makebox(0,0)[lb]{\smash{{{\SetFigFont{12}{14.4}{rm}.}}}}}
\put(455,125){\makebox(0,0)[lb]{\smash{{{\SetFigFont{12}{14.4}{rm}=}}}}}
\path(340,40)(340,60)(320,60)
	(320,40)(340,40)
\path(360,40)(360,60)(340,60)
	(340,40)(360,40)
\path(340,20)(340,40)(320,40)
	(320,20)(340,20)
\path(360,20)(360,40)(340,40)
	(340,20)(360,20)
\path(360,40)(360,40)(360,40)
	(360,40)(360,40)
\path(380,20)(380,40)(360,40)
	(360,20)(380,20)
\path(400,20)(400,40)(380,40)
	(380,20)(400,20)
\path(400,0)(400,20)(380,20)
	(380,0)(400,0)
\path(360,20)(360,20)(360,20)
	(360,20)(360,20)
\path(380,0)(380,20)(360,20)
	(360,0)(380,0)
\path(360,0)(360,20)(340,20)
	(340,0)(360,0)
\path(340,0)(340,20)(320,20)
	(320,0)(340,0)
\put(325,45){\makebox(0,0)[lb]{\smash{{{\SetFigFont{12}{14.4}{rm}3}}}}}
\put(345,45){\makebox(0,0)[lb]{\smash{{{\SetFigFont{12}{14.4}{rm}3}}}}}
\put(325,25){\makebox(0,0)[lb]{\smash{{{\SetFigFont{12}{14.4}{rm}2}}}}}
\put(345,25){\makebox(0,0)[lb]{\smash{{{\SetFigFont{12}{14.4}{rm}2}}}}}
\put(365,25){\makebox(0,0)[lb]{\smash{{{\SetFigFont{12}{14.4}{rm}2}}}}}
\put(385,25){\makebox(0,0)[lb]{\smash{{{\SetFigFont{12}{14.4}{rm}2}}}}}
\put(325,5){\makebox(0,0)[lb]{\smash{{{\SetFigFont{12}{14.4}{rm}1}}}}}
\put(345,5){\makebox(0,0)[lb]{\smash{{{\SetFigFont{12}{14.4}{rm}1}}}}}
\put(365,5){\makebox(0,0)[lb]{\smash{{{\SetFigFont{12}{14.4}{rm}1}}}}}
\put(385,5){\makebox(0,0)[lb]{\smash{{{\SetFigFont{12}{14.4}{rm}1}}}}}
\put(410,0){\makebox(0,0)[lb]{\smash{{{\SetFigFont{12}{14.4}{rm}.}}}}}
\path(170,20)(170,40)(150,40)
	(150,20)(170,20)
\path(190,20)(190,40)(170,40)
	(170,20)(190,20)
\path(190,0)(190,20)(170,20)
	(170,0)(190,0)
\path(210,0)(210,20)(190,20)
	(190,0)(210,0)
\path(230,0)(230,20)(210,20)
	(210,0)(230,0)
\path(270,20)(270,40)(250,40)
	(250,20)(270,20)
\path(270,0)(270,20)(250,20)
	(250,0)(270,0)
\path(290,20)(290,40)(270,40)
	(270,20)(290,20)
\path(290,0)(290,20)(270,20)
	(270,0)(290,0)
\path(170,0)(170,20)(150,20)
	(150,0)(170,0)
\put(155,25){\makebox(0,0)[lb]{\smash{{{\SetFigFont{12}{14.4}{rm}2}}}}}
\put(175,25){\makebox(0,0)[lb]{\smash{{{\SetFigFont{12}{14.4}{rm}2}}}}}
\put(155,5){\makebox(0,0)[lb]{\smash{{{\SetFigFont{12}{14.4}{rm}1}}}}}
\put(175,5){\makebox(0,0)[lb]{\smash{{{\SetFigFont{12}{14.4}{rm}1}}}}}
\put(195,5){\makebox(0,0)[lb]{\smash{{{\SetFigFont{12}{14.4}{rm}3}}}}}
\put(215,5){\makebox(0,0)[lb]{\smash{{{\SetFigFont{12}{14.4}{rm}3}}}}}
\put(255,25){\makebox(0,0)[lb]{\smash{{{\SetFigFont{12}{14.4}{rm}2}}}}}
\put(275,25){\makebox(0,0)[lb]{\smash{{{\SetFigFont{12}{14.4}{rm}2}}}}}
\put(255,5){\makebox(0,0)[lb]{\smash{{{\SetFigFont{12}{14.4}{rm}1}}}}}
\put(275,5){\makebox(0,0)[lb]{\smash{{{\SetFigFont{12}{14.4}{rm}1}}}}}
\put(240,5){\makebox(0,0)[lb]{\smash{{{\SetFigFont{12}{14.4}{rm}.}}}}}
\put(300,5){\makebox(0,0)[lb]{\smash{{{\SetFigFont{12}{14.4}{rm}=}}}}}
\end{picture}
\end{center}

\noindent
so that
$$
\K_{(4,2),(2,2,2)}(x_1,x_2)=x_1x_1 + x_2x_1 + x_2x_2 = s_{(2)}(x_1,x_2) \ .
$$
}
\end{example}

\Proof The proof relies on the fact that for a tableau $t$ of
rectangular weight $(k^{n+1})$, $d_i(t)$ is equal to the greatest
integer $s$ such that $\epsilon_i^s(t)\not =0$.
Suppose that the rightmost occurence of a letter $a_{i+1}$ in $t$
is {\it free}, {\it i.e.} can be transformed into $a_i$ by an application
of $\epsilon_i$. A tableau $u$ such that $t\,u$ is Yamanouchi will
then necessarily contain the column $c_i=a_i\cdots a_2a_1$.
Now the tableau $t\,c_i$ has one violation of the Yamanouchi condition
less than $t$, and $d_i(t\,c_i)=d_i(t)-1$. Iterating the process,
one arrives at a Yamanouchi tableau $t\,u$, and the minimal tableau $u$
constructed
is this way is clearly unique. The class of $u$ in the plactic monoid
is a product of mutually commuting columns $c_i$, and each $c_i$ appears
with multiplicity $d_i(t)$. \cqfd

\begin{example}{\rm With
\begin{center}
\setlength{\unitlength}{0.0085in}
\begingroup\makeatletter\ifx\SetFigFont\undefined
\def\x#1#2#3#4#5#6#7\relax{\def\x{#1#2#3#4#5#6}}%
\expandafter\x\fmtname xxxxxx\relax \def\y{splain}%
\ifx\x\y   
\gdef\SetFigFont#1#2#3{%
  \ifnum #1<17\tiny\else \ifnum #1<20\small\else
  \ifnum #1<24\normalsize\else \ifnum #1<29\large\else
  \ifnum #1<34\Large\else \ifnum #1<41\LARGE\else
     \huge\fi\fi\fi\fi\fi\fi
  \csname #3\endcsname}%
\else
\gdef\SetFigFont#1#2#3{\begingroup
  \count@#1\relax \ifnum 25<\count@\count@25\fi
  \def\x{\endgroup\@setsize\SetFigFont{#2pt}}%
  \expandafter\x
    \csname \romannumeral\the\count@ pt\expandafter\endcsname
    \csname @\romannumeral\the\count@ pt\endcsname
  \csname #3\endcsname}%
\fi
\fi\endgroup
\begin{picture}(135,57)(0,-10)
\path(55,20)(55,40)(35,40)
	(35,20)(55,20)
\path(55,0)(55,20)(35,20)
	(35,0)(55,0)
\path(75,0)(75,20)(55,20)
	(55,0)(75,0)
\path(95,0)(95,20)(75,20)
	(75,0)(95,0)
\path(115,0)(115,20)(95,20)
	(95,0)(115,0)
\path(135,0)(135,20)(115,20)
	(115,0)(135,0)
\put(40,25){\makebox(0,0)[lb]{\smash{{{\SetFigFont{12}{14.4}{rm}2}}}}}
\put(40,5){\makebox(0,0)[lb]{\smash{{{\SetFigFont{12}{14.4}{rm}1}}}}}
\put(60,5){\makebox(0,0)[lb]{\smash{{{\SetFigFont{12}{14.4}{rm}1}}}}}
\put(80,5){\makebox(0,0)[lb]{\smash{{{\SetFigFont{12}{14.4}{rm}2}}}}}
\put(100,5){\makebox(0,0)[lb]{\smash{{{\SetFigFont{12}{14.4}{rm}3}}}}}
\put(120,5){\makebox(0,0)[lb]{\smash{{{\SetFigFont{12}{14.4}{rm}3}}}}}
\put(0,5){\makebox(0,0)[lb]{\smash{{{\SetFigFont{12}{14.4}{rm}t}}}}}
\put(15,5){\makebox(0,0)[lb]{\smash{{{\SetFigFont{12}{14.4}{rm}=}}}}}
\end{picture}
\end{center}

\noindent
one has

\begin{center}
\setlength{\unitlength}{0.0085in}
\begingroup\makeatletter\ifx\SetFigFont\undefined
\def\x#1#2#3#4#5#6#7\relax{\def\x{#1#2#3#4#5#6}}%
\expandafter\x\fmtname xxxxxx\relax \def\y{splain}%
\ifx\x\y   
\gdef\SetFigFont#1#2#3{%
  \ifnum #1<17\tiny\else \ifnum #1<20\small\else
  \ifnum #1<24\normalsize\else \ifnum #1<29\large\else
  \ifnum #1<34\Large\else \ifnum #1<41\LARGE\else
     \huge\fi\fi\fi\fi\fi\fi
  \csname #3\endcsname}%
\else
\gdef\SetFigFont#1#2#3{\begingroup
  \count@#1\relax \ifnum 25<\count@\count@25\fi
  \def\x{\endgroup\@setsize\SetFigFont{#2pt}}%
  \expandafter\x
    \csname \romannumeral\the\count@ pt\expandafter\endcsname
    \csname @\romannumeral\the\count@ pt\endcsname
  \csname #3\endcsname}%
\fi
\fi\endgroup
\begin{picture}(320,77)(0,-10)
\path(20,20)(20,40)(0,40)
	(0,20)(20,20)
\path(20,0)(20,20)(0,20)
	(0,0)(20,0)
\path(40,0)(40,20)(20,20)
	(20,0)(40,0)
\path(60,0)(60,20)(40,20)
	(40,0)(60,0)
\path(80,0)(80,20)(60,20)
	(60,0)(80,0)
\path(100,0)(100,20)(80,20)
	(80,0)(100,0)
\path(140,20)(140,40)(120,40)
	(120,20)(140,20)
\path(160,20)(160,40)(140,40)
	(140,20)(160,20)
\path(140,0)(140,20)(120,20)
	(120,0)(140,0)
\path(160,0)(160,20)(140,20)
	(140,0)(160,0)
\path(180,0)(180,20)(160,20)
	(160,0)(180,0)
\put(5,25){\makebox(0,0)[lb]{\smash{{{\SetFigFont{12}{14.4}{rm}2}}}}}
\put(5,5){\makebox(0,0)[lb]{\smash{{{\SetFigFont{12}{14.4}{rm}1}}}}}
\put(25,5){\makebox(0,0)[lb]{\smash{{{\SetFigFont{12}{14.4}{rm}1}}}}}
\put(45,5){\makebox(0,0)[lb]{\smash{{{\SetFigFont{12}{14.4}{rm}2}}}}}
\put(65,5){\makebox(0,0)[lb]{\smash{{{\SetFigFont{12}{14.4}{rm}3}}}}}
\put(85,5){\makebox(0,0)[lb]{\smash{{{\SetFigFont{12}{14.4}{rm}3}}}}}
\path(240,40)(240,60)(220,60)
	(220,40)(240,40)
\path(260,40)(260,60)(240,60)
	(240,40)(260,40)
\path(240,20)(240,40)(220,40)
	(220,20)(240,20)
\path(260,20)(260,40)(240,40)
	(240,20)(260,20)
\path(280,20)(280,40)(260,40)
	(260,20)(280,20)
\path(240,0)(240,20)(220,20)
	(220,0)(240,0)
\path(260,0)(260,20)(240,20)
	(240,0)(260,0)
\path(280,0)(280,20)(260,20)
	(260,0)(280,0)
\path(300,0)(300,20)(280,20)
	(280,0)(300,0)
\path(320,0)(320,20)(300,20)
	(300,0)(320,0)
\path(300,20)(300,40)(280,40)
	(280,20)(300,20)
\put(110,5){\makebox(0,0)[lb]{\smash{{{\SetFigFont{12}{14.4}{rm}.}}}}}
\put(125,25){\makebox(0,0)[lb]{\smash{{{\SetFigFont{12}{14.4}{rm}2}}}}}
\put(145,25){\makebox(0,0)[lb]{\smash{{{\SetFigFont{12}{14.4}{rm}2}}}}}
\put(125,5){\makebox(0,0)[lb]{\smash{{{\SetFigFont{12}{14.4}{rm}1}}}}}
\put(145,5){\makebox(0,0)[lb]{\smash{{{\SetFigFont{12}{14.4}{rm}1}}}}}
\put(165,5){\makebox(0,0)[lb]{\smash{{{\SetFigFont{12}{14.4}{rm}1}}}}}
\put(195,5){\makebox(0,0)[lb]{\smash{{{\SetFigFont{12}{14.4}{rm}=}}}}}
\put(225,45){\makebox(0,0)[lb]{\smash{{{\SetFigFont{12}{14.4}{rm}3}}}}}
\put(245,45){\makebox(0,0)[lb]{\smash{{{\SetFigFont{12}{14.4}{rm}3}}}}}
\put(225,25){\makebox(0,0)[lb]{\smash{{{\SetFigFont{12}{14.4}{rm}2}}}}}
\put(245,25){\makebox(0,0)[lb]{\smash{{{\SetFigFont{12}{14.4}{rm}2}}}}}
\put(265,25){\makebox(0,0)[lb]{\smash{{{\SetFigFont{12}{14.4}{rm}2}}}}}
\put(285,25){\makebox(0,0)[lb]{\smash{{{\SetFigFont{12}{14.4}{rm}2}}}}}
\put(225,5){\makebox(0,0)[lb]{\smash{{{\SetFigFont{12}{14.4}{rm}1}}}}}
\put(245,5){\makebox(0,0)[lb]{\smash{{{\SetFigFont{12}{14.4}{rm}1}}}}}
\put(265,5){\makebox(0,0)[lb]{\smash{{{\SetFigFont{12}{14.4}{rm}1}}}}}
\put(285,5){\makebox(0,0)[lb]{\smash{{{\SetFigFont{12}{14.4}{rm}1}}}}}
\put(305,5){\makebox(0,0)[lb]{\smash{{{\SetFigFont{12}{14.4}{rm}1}}}}}
\end{picture}
\end{center}

\noindent
which corresponds to the following pairings

\begin{center}
\setlength{\unitlength}{0.0085in}
\begingroup\makeatletter\ifx\SetFigFont\undefined
\def\x#1#2#3#4#5#6#7\relax{\def\x{#1#2#3#4#5#6}}%
\expandafter\x\fmtname xxxxxx\relax \def\y{splain}%
\ifx\x\y   
\gdef\SetFigFont#1#2#3{%
  \ifnum #1<17\tiny\else \ifnum #1<20\small\else
  \ifnum #1<24\normalsize\else \ifnum #1<29\large\else
  \ifnum #1<34\Large\else \ifnum #1<41\LARGE\else
     \huge\fi\fi\fi\fi\fi\fi
  \csname #3\endcsname}%
\else
\gdef\SetFigFont#1#2#3{\begingroup
  \count@#1\relax \ifnum 25<\count@\count@25\fi
  \def\x{\endgroup\@setsize\SetFigFont{#2pt}}%
  \expandafter\x
    \csname \romannumeral\the\count@ pt\expandafter\endcsname
    \csname @\romannumeral\the\count@ pt\endcsname
  \csname #3\endcsname}%
\fi
\fi\endgroup
\begin{picture}(227,65)(0,-10)
\put(0,25){\makebox(0,0)[lb]{\smash{{{\SetFigFont{12}{14.4}{rm}2}}}}}
\put(20,25){\makebox(0,0)[lb]{\smash{{{\SetFigFont{12}{14.4}{rm}1}}}}}
\put(40,25){\makebox(0,0)[lb]{\smash{{{\SetFigFont{12}{14.4}{rm}1}}}}}
\put(60,25){\makebox(0,0)[lb]{\smash{{{\SetFigFont{12}{14.4}{rm}2}}}}}
\put(80,25){\makebox(0,0)[lb]{\smash{{{\SetFigFont{12}{14.4}{rm}3}}}}}
\put(100,25){\makebox(0,0)[lb]{\smash{{{\SetFigFont{12}{14.4}{rm}3}}}}}
\put(125,25){\makebox(0,0)[lb]{\smash{{{\SetFigFont{12}{14.4}{rm}.}}}}}
\put(140,25){\makebox(0,0)[lb]{\smash{{{\SetFigFont{12}{14.4}{rm}2}}}}}
\put(160,25){\makebox(0,0)[lb]{\smash{{{\SetFigFont{12}{14.4}{rm}1}}}}}
\put(180,25){\makebox(0,0)[lb]{\smash{{{\SetFigFont{12}{14.4}{rm}2}}}}}
\put(200,25){\makebox(0,0)[lb]{\smash{{{\SetFigFont{12}{14.4}{rm}1}}}}}
\put(220,25){\makebox(0,0)[lb]{\smash{{{\SetFigFont{12}{14.4}{rm}1}}}}}
\path(145,20)(145,10)(160,10)
	(165,10)(165,20)
\path(105,20)(105,10)(140,10)(140,20)
\path(85,40)(85,50)(180,50)(180,40)
\path(5,40)(5,50)(25,50)(25,40)
\path(185,40)(185,50)(205,50)(205,40)
\path(65,20)(65,0)(225,0)(225,20)
\end{picture}
\end{center}

\noindent
The monomial associated to $t$ is $x_2x_2x_1$, obtained by reading
the lengths of the columns of $u$.
}
\end{example}

If $t'=\Omega_2(t)$, one has $d_i(t') = d_{n-i+1}(t)$, so that
the polynomials $\K_{\lambda \mu}$ are symmetrical in any pair of variables
with
complementary indices $(x_i,x_{n+1-i})$. For certain families of
partitions, they are even symmetric in the whole set of variables
$x_1,x_2,\ldots,x_n$. It will be convenient to adopt the
following notations:

Let $k$ be an integer ${\ge} 1$, $\beta$ a partition of length
$r\le n$ whose parts are $\le k$, and $\alpha$ a partition
of the same weight than $\beta$, and length $\le n+1$. Define a
partition $\lambda=[\alpha,\beta]_{n+1}^k$ of length $n+1$ by
$$
\lambda=(\alpha_1+k,\ldots,\alpha_r+k,k,k,\ldots,k,k-\beta_s,\ldots,k-\beta_1)\
{}.
$$
When $k=\beta_1$ we simply write $\lambda=[\alpha,\beta]_{n+1}$.

\begin{theorem}\label{RECT}
{\rm (i)} Using these notations, when the partition
$\alpha$ is a row $\alpha=(m)$, we have
$$\K_{[(m),\beta]_{n+1}^k,(k^{n+1})}(x_1,\ldots,x_n)=
s_\beta(x_1,\ldots,x_n)\,.$$

{\rm (ii)} Similarly, when $\beta=(m)$,
$\K_{[\alpha,(m)]_{n+1}^k,(k^{n+1})}(x_1,\ldots,x_n)
=s_\alpha(x_1,\ldots,x_n)\,.$

{\rm (iii)} In any case, for $k \ge \max(\alpha_1,\beta_1)$,
$$\K_{[\alpha,\beta]_{n+1}^k,(k^{n+1})}(x_1,\ldots,x_n)=
\K_{[\beta,\alpha]_{n+1}^k,(k^{n+1})}(x_1,\ldots,x_n)\,.$$
\end{theorem}

\Proof Set $\lambda = [\alpha,\beta]_{n+1}^k$.
Since one clearly has
\begin{equation}\label{STAB}
\K_{\lambda+(r^{n+1}),((k+r)^{n+1})} = \K_{\lambda,(k^{n+1})}
\end{equation}
one may always assume that $k \ge \max(\alpha_1,\beta_1)$, and (i) will
be a consequence of (ii) and (iii).

To establish (ii), one remarks that
when $\beta=(m)$ is a row, the free letters of a tableau
$t\in\tab(\lambda,(k^{n+1}))$ are exactly those filling the piece of
shape $\alpha$ lying outside the rectangle, and this piece
can be any tableau of weight $m$ over ${2,3,\ldots,n+1}$.
Hence $\K_{[\alpha,(m)]_{n+1}^k,(k^{n+1})}(x_1,\ldots,x_n)
=s_\alpha(x_1,\ldots,x_n)\,.$

Finally, point (iii) follows from the existence of an exponent-preserving
map on tableaux exchanging
$\tab([\alpha,\beta]^k_{n+1},(k^{n+1}))$ and
$\tab([\beta,\alpha]^k_{n+1},(k^{n+1}))$.

Let $t\in\tab([\alpha,\beta]^k_{n+1},(k^{n+1}))$. This map is defined for
$k$ sufficiently large ($k\ge\max(\alpha_1,\beta_1))$, but one can also
consider it as a well-defined involution on $\sl_{n+1}$-tableaux (which
are ordinary tableaux modulo the relation $a_{n+1}a_n\cdots a_1\equiv 1$).
To construct the
image $t'$  of $t$ by this involution, complete the part of $t$ lying
in the rectangle $(k^{n+1})$ by writing on the top of each column
its complementary subset, in decreasing order from bottom to top.
One obtains in this way the mirror image of a tableau of shape
$\beta$. Put this tableau at the right of an initially empty
rectangle $(k^{n+1})$, and fill its top right corner by
the mirror image of the part of $t$ lying inside $\alpha$.
Now fill the columns by the complementary
subsets of the corresponding columns, and finally erase the mirror
$\alpha$-tableau
of the top right corner. This correspondence is obviously self-inverse
(for $k$ large enough), and it is not difficult to check that
it preserves the exponents. For example,
with $k=5$, $n=4$, $\alpha=(4,2,2)$, $\beta=(5,3)$ and
\begin{center}
\setlength{\unitlength}{0.0085in}
\begingroup\makeatletter\ifx\SetFigFont\undefined
\def\x#1#2#3#4#5#6#7\relax{\def\x{#1#2#3#4#5#6}}%
\expandafter\x\fmtname xxxxxx\relax \def\y{splain}%
\ifx\x\y   
\gdef\SetFigFont#1#2#3{%
  \ifnum #1<17\tiny\else \ifnum #1<20\small\else
  \ifnum #1<24\normalsize\else \ifnum #1<29\large\else
  \ifnum #1<34\Large\else \ifnum #1<41\LARGE\else
     \huge\fi\fi\fi\fi\fi\fi
  \csname #3\endcsname}%
\else
\gdef\SetFigFont#1#2#3{\begingroup
  \count@#1\relax \ifnum 25<\count@\count@25\fi
  \def\x{\endgroup\@setsize\SetFigFont{#2pt}}%
  \expandafter\x
    \csname \romannumeral\the\count@ pt\expandafter\endcsname
    \csname @\romannumeral\the\count@ pt\endcsname
  \csname #3\endcsname}%
\fi
\fi\endgroup
\begin{picture}(215,97)(0,-10)
\path(55,60)(55,80)(35,80)
	(35,60)(55,60)
\path(75,60)(75,80)(55,80)
	(55,60)(75,60)
\path(55,40)(55,60)(35,60)
	(35,40)(55,40)
\path(75,40)(75,60)(55,60)
	(55,40)(75,40)
\path(95,40)(95,60)(75,60)
	(75,40)(95,40)
\path(115,40)(115,60)(95,60)
	(95,40)(115,40)
\path(135,40)(135,60)(115,60)
	(115,40)(135,40)
\path(155,40)(155,60)(135,60)
	(135,40)(155,40)
\path(175,40)(175,60)(155,60)
	(155,40)(175,40)
\path(55,20)(55,40)(35,40)
	(35,20)(55,20)
\path(75,20)(75,40)(55,40)
	(55,20)(75,20)
\path(95,20)(95,40)(75,40)
	(75,20)(95,20)
\path(95,20)(95,40)(95,40)
	(95,20)(95,20)
\path(115,20)(115,20)(115,20)
	(115,20)(115,20)
\path(115,20)(115,40)(95,40)
	(95,20)(115,20)
\path(135,20)(135,40)(115,40)
	(115,20)(135,20)
\path(155,20)(155,40)(135,40)
	(135,20)(155,20)
\path(175,20)(175,40)(155,40)
	(155,20)(175,20)
\path(55,0)(55,20)(35,20)
	(35,0)(55,0)
\path(75,0)(75,20)(55,20)
	(55,0)(75,0)
\path(95,0)(95,20)(75,20)
	(75,0)(95,0)
\path(115,0)(115,20)(95,20)
	(95,0)(115,0)
\path(135,0)(135,20)(115,20)
	(115,0)(135,0)
\path(155,0)(155,20)(135,20)
	(135,0)(155,0)
\path(175,0)(175,20)(155,20)
	(155,0)(175,0)
\path(195,0)(195,20)(175,20)
	(175,0)(195,0)
\path(215,0)(215,20)(195,20)
	(195,0)(215,0)
\put(0,5){\makebox(0,0)[lb]{\smash{{{\SetFigFont{12}{14.4}{rm}t  =}}}}}
\put(40,65){\makebox(0,0)[lb]{\smash{{{\SetFigFont{12}{14.4}{rm}5}}}}}
\put(60,65){\makebox(0,0)[lb]{\smash{{{\SetFigFont{12}{14.4}{rm}5}}}}}
\put(40,45){\makebox(0,0)[lb]{\smash{{{\SetFigFont{12}{14.4}{rm}3}}}}}
\put(60,45){\makebox(0,0)[lb]{\smash{{{\SetFigFont{12}{14.4}{rm}3}}}}}
\put(80,45){\makebox(0,0)[lb]{\smash{{{\SetFigFont{12}{14.4}{rm}4}}}}}
\put(100,45){\makebox(0,0)[lb]{\smash{{{\SetFigFont{12}{14.4}{rm}4}}}}}
\put(120,45){\makebox(0,0)[lb]{\smash{{{\SetFigFont{12}{14.4}{rm}4}}}}}
\put(140,45){\makebox(0,0)[lb]{\smash{{{\SetFigFont{12}{14.4}{bf}5}}}}}
\put(160,45){\makebox(0,0)[lb]{\smash{{{\SetFigFont{12}{14.4}{bf}5}}}}}
\put(40,25){\makebox(0,0)[lb]{\smash{{{\SetFigFont{12}{14.4}{rm}2}}}}}
\put(60,25){\makebox(0,0)[lb]{\smash{{{\SetFigFont{12}{14.4}{rm}2}}}}}
\put(80,25){\makebox(0,0)[lb]{\smash{{{\SetFigFont{12}{14.4}{rm}2}}}}}
\put(100,25){\makebox(0,0)[lb]{\smash{{{\SetFigFont{12}{14.4}{rm}3}}}}}
\put(120,25){\makebox(0,0)[lb]{\smash{{{\SetFigFont{12}{14.4}{rm}3}}}}}
\put(40,5){\makebox(0,0)[lb]{\smash{{{\SetFigFont{12}{14.4}{rm}1}}}}}
\put(60,5){\makebox(0,0)[lb]{\smash{{{\SetFigFont{12}{14.4}{rm}1}}}}}
\put(80,5){\makebox(0,0)[lb]{\smash{{{\SetFigFont{12}{14.4}{rm}1}}}}}
\put(100,5){\makebox(0,0)[lb]{\smash{{{\SetFigFont{12}{14.4}{rm}1}}}}}
\put(120,5){\makebox(0,0)[lb]{\smash{{{\SetFigFont{12}{14.4}{rm}1}}}}}
\put(140,25){\makebox(0,0)[lb]{\smash{{{\SetFigFont{12}{14.4}{bf}3}}}}}
\put(160,25){\makebox(0,0)[lb]{\smash{{{\SetFigFont{12}{14.4}{bf}4}}}}}
\put(140,5){\makebox(0,0)[lb]{\smash{{{\SetFigFont{12}{14.4}{bf}2}}}}}
\put(160,5){\makebox(0,0)[lb]{\smash{{{\SetFigFont{12}{14.4}{bf}2}}}}}
\put(180,5){\makebox(0,0)[lb]{\smash{{{\SetFigFont{12}{14.4}{bf}4}}}}}
\put(200,5){\makebox(0,0)[lb]{\smash{{{\SetFigFont{12}{14.4}{bf}5}}}}}
\end{picture}
\end{center}
one constructs sucessively
\begin{center}
\setlength{\unitlength}{0.0085in}
\begingroup\makeatletter\ifx\SetFigFont\undefined
\def\x#1#2#3#4#5#6#7\relax{\def\x{#1#2#3#4#5#6}}%
\expandafter\x\fmtname xxxxxx\relax \def\y{splain}%
\ifx\x\y   
\gdef\SetFigFont#1#2#3{%
  \ifnum #1<17\tiny\else \ifnum #1<20\small\else
  \ifnum #1<24\normalsize\else \ifnum #1<29\large\else
  \ifnum #1<34\Large\else \ifnum #1<41\LARGE\else
     \huge\fi\fi\fi\fi\fi\fi
  \csname #3\endcsname}%
\else
\gdef\SetFigFont#1#2#3{\begingroup
  \count@#1\relax \ifnum 25<\count@\count@25\fi
  \def\x{\endgroup\@setsize\SetFigFont{#2pt}}%
  \expandafter\x
    \csname \romannumeral\the\count@ pt\expandafter\endcsname
    \csname @\romannumeral\the\count@ pt\endcsname
  \csname #3\endcsname}%
\fi
\fi\endgroup
\begin{picture}(480,117)(0,-10)
\path(20,60)(20,80)(0,80)
	(0,60)(20,60)
\path(40,60)(40,80)(20,80)
	(20,60)(40,60)
\path(20,40)(20,60)(0,60)
	(0,40)(20,40)
\path(40,40)(40,60)(20,60)
	(20,40)(40,40)
\path(60,40)(60,60)(40,60)
	(40,40)(60,40)
\path(80,40)(80,60)(60,60)
	(60,40)(80,40)
\path(100,40)(100,60)(80,60)
	(80,40)(100,40)
\path(120,40)(120,60)(100,60)
	(100,40)(120,40)
\path(140,40)(140,60)(120,60)
	(120,40)(140,40)
\path(20,20)(20,40)(0,40)
	(0,20)(20,20)
\path(40,20)(40,40)(20,40)
	(20,20)(40,20)
\path(60,20)(60,40)(40,40)
	(40,20)(60,20)
\path(60,20)(60,40)(60,40)
	(60,20)(60,20)
\path(80,20)(80,20)(80,20)
	(80,20)(80,20)
\path(80,20)(80,40)(60,40)
	(60,20)(80,20)
\path(100,20)(100,40)(80,40)
	(80,20)(100,20)
\path(120,20)(120,40)(100,40)
	(100,20)(120,20)
\path(140,20)(140,40)(120,40)
	(120,20)(140,20)
\path(20,0)(20,20)(0,20)
	(0,0)(20,0)
\path(40,0)(40,20)(20,20)
	(20,0)(40,0)
\path(60,0)(60,20)(40,20)
	(40,0)(60,0)
\path(80,0)(80,20)(60,20)
	(60,0)(80,0)
\path(100,0)(100,20)(80,20)
	(80,0)(100,0)
\path(120,0)(120,20)(100,20)
	(100,0)(120,0)
\path(140,0)(140,20)(120,20)
	(120,0)(140,0)
\path(160,0)(160,20)(140,20)
	(140,0)(160,0)
\path(180,0)(180,20)(160,20)
	(160,0)(180,0)
\dashline{4.000}(0,85)(0,100)(100,100)(100,65)
\path(300,0)(300,20)(280,20)
	(280,0)(300,0)
\path(320,0)(320,20)(300,20)
	(300,0)(320,0)
\path(340,0)(340,20)(320,20)
	(320,0)(340,0)
\path(360,0)(360,20)(340,20)
	(340,0)(360,0)
\path(380,0)(380,20)(360,20)
	(360,0)(380,0)
\path(400,0)(400,20)(380,20)
	(380,0)(400,0)
\path(420,0)(420,20)(400,20)
	(400,0)(420,0)
\path(440,0)(440,20)(420,20)
	(420,0)(440,0)
\path(460,0)(460,20)(440,20)
	(440,0)(460,0)
\path(480,0)(480,20)(460,20)
	(460,0)(480,0)
\path(300,20)(300,40)(280,40)
	(280,20)(300,20)
\path(320,20)(320,40)(300,40)
	(300,20)(320,20)
\path(340,20)(340,40)(320,40)
	(320,20)(340,20)
\path(360,20)(360,40)(340,40)
	(340,20)(360,20)
\path(380,20)(380,40)(360,40)
	(360,20)(380,20)
\path(400,20)(400,40)(380,40)
	(380,20)(400,20)
\path(420,20)(420,40)(400,40)
	(400,20)(420,20)
\path(440,20)(440,40)(420,40)
	(420,20)(440,20)
\path(300,40)(300,60)(280,60)
	(280,40)(300,40)
\path(320,40)(320,60)(300,60)
	(300,40)(320,40)
\path(340,40)(340,60)(320,60)
	(320,40)(340,40)
\path(300,60)(300,80)(280,80)
	(280,60)(300,60)
\path(320,60)(320,80)(300,80)
	(300,60)(320,60)
\path(340,60)(340,80)(320,80)
	(320,60)(340,60)
\path(300,80)(300,100)(280,100)
	(280,80)(300,80)
\dashline{4.000}(300,100)(380,100)(380,40)
\put(5,65){\makebox(0,0)[lb]{\smash{{{\SetFigFont{12}{14.4}{rm}5}}}}}
\put(25,65){\makebox(0,0)[lb]{\smash{{{\SetFigFont{12}{14.4}{rm}5}}}}}
\put(5,45){\makebox(0,0)[lb]{\smash{{{\SetFigFont{12}{14.4}{rm}3}}}}}
\put(25,45){\makebox(0,0)[lb]{\smash{{{\SetFigFont{12}{14.4}{rm}3}}}}}
\put(45,45){\makebox(0,0)[lb]{\smash{{{\SetFigFont{12}{14.4}{rm}4}}}}}
\put(65,45){\makebox(0,0)[lb]{\smash{{{\SetFigFont{12}{14.4}{rm}4}}}}}
\put(85,45){\makebox(0,0)[lb]{\smash{{{\SetFigFont{12}{14.4}{rm}4}}}}}
\put(105,45){\makebox(0,0)[lb]{\smash{{{\SetFigFont{12}{14.4}{rm}5}}}}}
\put(125,45){\makebox(0,0)[lb]{\smash{{{\SetFigFont{12}{14.4}{rm}5}}}}}
\put(5,25){\makebox(0,0)[lb]{\smash{{{\SetFigFont{12}{14.4}{rm}2}}}}}
\put(25,25){\makebox(0,0)[lb]{\smash{{{\SetFigFont{12}{14.4}{rm}2}}}}}
\put(45,25){\makebox(0,0)[lb]{\smash{{{\SetFigFont{12}{14.4}{rm}2}}}}}
\put(65,25){\makebox(0,0)[lb]{\smash{{{\SetFigFont{12}{14.4}{rm}3}}}}}
\put(85,25){\makebox(0,0)[lb]{\smash{{{\SetFigFont{12}{14.4}{rm}3}}}}}
\put(5,5){\makebox(0,0)[lb]{\smash{{{\SetFigFont{12}{14.4}{rm}1}}}}}
\put(25,5){\makebox(0,0)[lb]{\smash{{{\SetFigFont{12}{14.4}{rm}1}}}}}
\put(45,5){\makebox(0,0)[lb]{\smash{{{\SetFigFont{12}{14.4}{rm}1}}}}}
\put(65,5){\makebox(0,0)[lb]{\smash{{{\SetFigFont{12}{14.4}{rm}1}}}}}
\put(85,5){\makebox(0,0)[lb]{\smash{{{\SetFigFont{12}{14.4}{rm}1}}}}}
\put(105,25){\makebox(0,0)[lb]{\smash{{{\SetFigFont{12}{14.4}{rm}3}}}}}
\put(125,25){\makebox(0,0)[lb]{\smash{{{\SetFigFont{12}{14.4}{rm}4}}}}}
\put(105,5){\makebox(0,0)[lb]{\smash{{{\SetFigFont{12}{14.4}{rm}2}}}}}
\put(125,5){\makebox(0,0)[lb]{\smash{{{\SetFigFont{12}{14.4}{rm}2}}}}}
\put(145,5){\makebox(0,0)[lb]{\smash{{{\SetFigFont{12}{14.4}{rm}4}}}}}
\put(165,5){\makebox(0,0)[lb]{\smash{{{\SetFigFont{12}{14.4}{rm}5}}}}}
\put(5,85){\makebox(0,0)[lb]{\smash{{{\SetFigFont{12}{14.4}{bf}4}}}}}
\put(25,85){\makebox(0,0)[lb]{\smash{{{\SetFigFont{12}{14.4}{bf}4}}}}}
\put(45,85){\makebox(0,0)[lb]{\smash{{{\SetFigFont{12}{14.4}{bf}3}}}}}
\put(45,65){\makebox(0,0)[lb]{\smash{{{\SetFigFont{12}{14.4}{bf}5}}}}}
\put(65,85){\makebox(0,0)[lb]{\smash{{{\SetFigFont{12}{14.4}{bf}2}}}}}
\put(65,65){\makebox(0,0)[lb]{\smash{{{\SetFigFont{12}{14.4}{bf}5}}}}}
\put(85,85){\makebox(0,0)[lb]{\smash{{{\SetFigFont{12}{14.4}{bf}2}}}}}
\put(85,65){\makebox(0,0)[lb]{\smash{{{\SetFigFont{12}{14.4}{bf}5}}}}}
\put(305,85){\makebox(0,0)[lb]{\smash{{{\SetFigFont{12}{14.4}{bf}5}}}}}
\put(325,85){\makebox(0,0)[lb]{\smash{{{\SetFigFont{12}{14.4}{bf}4}}}}}
\put(345,85){\makebox(0,0)[lb]{\smash{{{\SetFigFont{12}{14.4}{bf}2}}}}}
\put(365,85){\makebox(0,0)[lb]{\smash{{{\SetFigFont{12}{14.4}{bf}2}}}}}
\put(345,65){\makebox(0,0)[lb]{\smash{{{\SetFigFont{12}{14.4}{bf}4}}}}}
\put(365,65){\makebox(0,0)[lb]{\smash{{{\SetFigFont{12}{14.4}{bf}3}}}}}
\put(345,45){\makebox(0,0)[lb]{\smash{{{\SetFigFont{12}{14.4}{bf}5}}}}}
\put(365,45){\makebox(0,0)[lb]{\smash{{{\SetFigFont{12}{14.4}{bf}5}}}}}
\put(285,85){\makebox(0,0)[lb]{\smash{{{\SetFigFont{12}{14.4}{rm}5}}}}}
\put(285,65){\makebox(0,0)[lb]{\smash{{{\SetFigFont{12}{14.4}{rm}4}}}}}
\put(305,65){\makebox(0,0)[lb]{\smash{{{\SetFigFont{12}{14.4}{rm}4}}}}}
\put(325,65){\makebox(0,0)[lb]{\smash{{{\SetFigFont{12}{14.4}{rm}5}}}}}
\put(285,45){\makebox(0,0)[lb]{\smash{{{\SetFigFont{12}{14.4}{rm}3}}}}}
\put(305,45){\makebox(0,0)[lb]{\smash{{{\SetFigFont{12}{14.4}{rm}3}}}}}
\put(325,45){\makebox(0,0)[lb]{\smash{{{\SetFigFont{12}{14.4}{rm}3}}}}}
\put(285,25){\makebox(0,0)[lb]{\smash{{{\SetFigFont{12}{14.4}{rm}2}}}}}
\put(305,25){\makebox(0,0)[lb]{\smash{{{\SetFigFont{12}{14.4}{rm}2}}}}}
\put(325,25){\makebox(0,0)[lb]{\smash{{{\SetFigFont{12}{14.4}{rm}2}}}}}
\put(345,25){\makebox(0,0)[lb]{\smash{{{\SetFigFont{12}{14.4}{rm}3}}}}}
\put(365,25){\makebox(0,0)[lb]{\smash{{{\SetFigFont{12}{14.4}{rm}4}}}}}
\put(385,25){\makebox(0,0)[lb]{\smash{{{\SetFigFont{12}{14.4}{rm}5}}}}}
\put(405,25){\makebox(0,0)[lb]{\smash{{{\SetFigFont{12}{14.4}{rm}5}}}}}
\put(425,25){\makebox(0,0)[lb]{\smash{{{\SetFigFont{12}{14.4}{rm}5}}}}}
\put(285,5){\makebox(0,0)[lb]{\smash{{{\SetFigFont{12}{14.4}{rm}1}}}}}
\put(305,5){\makebox(0,0)[lb]{\smash{{{\SetFigFont{12}{14.4}{rm}1}}}}}
\put(325,5){\makebox(0,0)[lb]{\smash{{{\SetFigFont{12}{14.4}{rm}1}}}}}
\put(345,5){\makebox(0,0)[lb]{\smash{{{\SetFigFont{12}{14.4}{rm}1}}}}}
\put(365,5){\makebox(0,0)[lb]{\smash{{{\SetFigFont{12}{14.4}{rm}1}}}}}
\put(385,5){\makebox(0,0)[lb]{\smash{{{\SetFigFont{12}{14.4}{rm}2}}}}}
\put(405,5){\makebox(0,0)[lb]{\smash{{{\SetFigFont{12}{14.4}{rm}2}}}}}
\put(425,5){\makebox(0,0)[lb]{\smash{{{\SetFigFont{12}{14.4}{rm}3}}}}}
\put(445,5){\makebox(0,0)[lb]{\smash{{{\SetFigFont{12}{14.4}{rm}4}}}}}
\put(465,5){\makebox(0,0)[lb]{\smash{{{\SetFigFont{12}{14.4}{rm}4}}}}}
\put(215,40){\makebox(0,0)[lb]{\smash{{{\SetFigFont{12}{14.4}{rm}and}}}}}
\end{picture}
\end{center}
so that
\begin{center}
\setlength{\unitlength}{0.0085in}
\begingroup\makeatletter\ifx\SetFigFont\undefined
\def\x#1#2#3#4#5#6#7\relax{\def\x{#1#2#3#4#5#6}}%
\expandafter\x\fmtname xxxxxx\relax \def\y{splain}%
\ifx\x\y   
\gdef\SetFigFont#1#2#3{%
  \ifnum #1<17\tiny\else \ifnum #1<20\small\else
  \ifnum #1<24\normalsize\else \ifnum #1<29\large\else
  \ifnum #1<34\Large\else \ifnum #1<41\LARGE\else
     \huge\fi\fi\fi\fi\fi\fi
  \csname #3\endcsname}%
\else
\gdef\SetFigFont#1#2#3{\begingroup
  \count@#1\relax \ifnum 25<\count@\count@25\fi
  \def\x{\endgroup\@setsize\SetFigFont{#2pt}}%
  \expandafter\x
    \csname \romannumeral\the\count@ pt\expandafter\endcsname
    \csname @\romannumeral\the\count@ pt\endcsname
  \csname #3\endcsname}%
\fi
\fi\endgroup
\begin{picture}(240,117)(0,-10)
\path(60,0)(60,20)(40,20)
	(40,0)(60,0)
\path(80,0)(80,20)(60,20)
	(60,0)(80,0)
\path(100,0)(100,20)(80,20)
	(80,0)(100,0)
\path(120,0)(120,20)(100,20)
	(100,0)(120,0)
\path(140,0)(140,20)(120,20)
	(120,0)(140,0)
\path(160,0)(160,20)(140,20)
	(140,0)(160,0)
\path(180,0)(180,20)(160,20)
	(160,0)(180,0)
\path(200,0)(200,20)(180,20)
	(180,0)(200,0)
\path(220,0)(220,20)(200,20)
	(200,0)(220,0)
\path(240,0)(240,20)(220,20)
	(220,0)(240,0)
\path(60,20)(60,40)(40,40)
	(40,20)(60,20)
\path(80,20)(80,40)(60,40)
	(60,20)(80,20)
\path(100,20)(100,40)(80,40)
	(80,20)(100,20)
\path(120,20)(120,40)(100,40)
	(100,20)(120,20)
\path(140,20)(140,40)(120,40)
	(120,20)(140,20)
\path(160,20)(160,40)(140,40)
	(140,20)(160,20)
\path(180,20)(180,40)(160,40)
	(160,20)(180,20)
\path(200,20)(200,40)(180,40)
	(180,20)(200,20)
\path(60,40)(60,60)(40,60)
	(40,40)(60,40)
\path(80,40)(80,60)(60,60)
	(60,40)(80,40)
\path(100,40)(100,60)(80,60)
	(80,40)(100,40)
\path(60,60)(60,80)(40,80)
	(40,60)(60,60)
\path(80,60)(80,80)(60,80)
	(60,60)(80,60)
\path(100,60)(100,80)(80,80)
	(80,60)(100,60)
\path(60,80)(60,100)(40,100)
	(40,80)(60,80)
\put(45,85){\makebox(0,0)[lb]{\smash{{{\SetFigFont{12}{14.4}{rm}5}}}}}
\put(45,65){\makebox(0,0)[lb]{\smash{{{\SetFigFont{12}{14.4}{rm}4}}}}}
\put(65,65){\makebox(0,0)[lb]{\smash{{{\SetFigFont{12}{14.4}{rm}4}}}}}
\put(85,65){\makebox(0,0)[lb]{\smash{{{\SetFigFont{12}{14.4}{rm}5}}}}}
\put(45,45){\makebox(0,0)[lb]{\smash{{{\SetFigFont{12}{14.4}{rm}3}}}}}
\put(65,45){\makebox(0,0)[lb]{\smash{{{\SetFigFont{12}{14.4}{rm}3}}}}}
\put(85,45){\makebox(0,0)[lb]{\smash{{{\SetFigFont{12}{14.4}{rm}3}}}}}
\put(45,25){\makebox(0,0)[lb]{\smash{{{\SetFigFont{12}{14.4}{rm}2}}}}}
\put(65,25){\makebox(0,0)[lb]{\smash{{{\SetFigFont{12}{14.4}{rm}2}}}}}
\put(85,25){\makebox(0,0)[lb]{\smash{{{\SetFigFont{12}{14.4}{rm}2}}}}}
\put(105,25){\makebox(0,0)[lb]{\smash{{{\SetFigFont{12}{14.4}{rm}3}}}}}
\put(125,25){\makebox(0,0)[lb]{\smash{{{\SetFigFont{12}{14.4}{rm}4}}}}}
\put(145,25){\makebox(0,0)[lb]{\smash{{{\SetFigFont{12}{14.4}{bf}5}}}}}
\put(165,25){\makebox(0,0)[lb]{\smash{{{\SetFigFont{12}{14.4}{bf}5}}}}}
\put(45,5){\makebox(0,0)[lb]{\smash{{{\SetFigFont{12}{14.4}{rm}1}}}}}
\put(65,5){\makebox(0,0)[lb]{\smash{{{\SetFigFont{12}{14.4}{rm}1}}}}}
\put(85,5){\makebox(0,0)[lb]{\smash{{{\SetFigFont{12}{14.4}{rm}1}}}}}
\put(105,5){\makebox(0,0)[lb]{\smash{{{\SetFigFont{12}{14.4}{rm}1}}}}}
\put(125,5){\makebox(0,0)[lb]{\smash{{{\SetFigFont{12}{14.4}{rm}1}}}}}
\put(185,25){\makebox(0,0)[lb]{\smash{{{\SetFigFont{12}{14.4}{bf}5}}}}}
\put(145,5){\makebox(0,0)[lb]{\smash{{{\SetFigFont{12}{14.4}{bf}2}}}}}
\put(165,5){\makebox(0,0)[lb]{\smash{{{\SetFigFont{12}{14.4}{bf}2}}}}}
\put(185,5){\makebox(0,0)[lb]{\smash{{{\SetFigFont{12}{14.4}{bf}3}}}}}
\put(205,5){\makebox(0,0)[lb]{\smash{{{\SetFigFont{12}{14.4}{bf}4}}}}}
\put(225,5){\makebox(0,0)[lb]{\smash{{{\SetFigFont{12}{14.4}{bf}4}}}}}
\put(0,40){\makebox(0,0)[lb]{\smash{{{\SetFigFont{12}{14.4}{rm}t'  =}}}}}
\end{picture}
\end{center}
One can check that both $t$ and $t'$ have exponents
${\bf d}(t)={\bf d}(t')=(2,\,1,\,2,\,1)$. \cqfd

One can also give an explicit formula in the case of the standard weight
$\mu=(1^{n+1})$.
Set $t_i=x_i\,(1-x_i)^{-1}$, $i{=}1,\,\ldots ,\,n$. Then:
\begin{equation}
{\K_{\lambda,(1^{n+1})}(x_1,\ldots,x_n) \over \prod_{1\le i \le n} (1-x_i)}
=
\sum_{\mu} K_{\lambda'\mu} \sum_{I\in \S(\mu)} t_{i_1}t_{i_1+i_2}\cdots
t_{i_1+\cdots +i_{r-1}} \,,
\end{equation}
where $\S(\mu)$ denotes the set of all compositions $I=(i_1,\ldots ,i_r)$
obtained
by permutation of the parts of $\mu=(\mu_1,\ldots ,\mu_r)$.

This formula is in fact the commutative image of an identity
in the algebra of noncommutative symmetric functions, defined in
\cite{GKLLRT}. Using the notations of \cite{GKLLRT}, one can
check that the generating function
$$
\sum_{\lambda}\K_{\lambda,(1^{n+1})}(x_1,\ldots,x_n)s_\lambda
$$
is the commutative image of
$$
{\cal K}_{n+1}(x_1,\ldots,x_n) =
\sum_{|I|=n+1} \left( \prod_{d\in D(I)} x_d \right) R_I
$$
where for a composition $I=(i_1,\ldots,i_r)$,
$D(I)=\{i_1,i_1+i_2,\ldots,i_1+\cdots+i_{r-1}\}$ and $R_I$
is the associated noncommutative ribbon Schur function.
Expanding on the basis of products of complete symmetric functions,
one finds
$$
{\cal K}_{n+1}=
\sum_{|J|=n+1} \left(
\sum_{I\ge J} (-1)^{\ell(J)-\ell(I)}
\prod_{d\in D(I)} x_d  \right)  S^J
$$
and it remains to show that the coefficient of $S^J$ is this expression
is equal to
$$
\prod_{d\in D(J)} x_d \prod_{e\not\in D(J)}(1-x_e) \ ,
$$
which follows from a straightforward argument of inclusion-exclusion.

\bigskip\noindent
{\bf Remarks}

\medskip\noindent
{\bf 1.} The stability property (\ref{STAB})
simply means that the $\K$-polynomials are associated to $\sl_{n+1}$
rather than to $\gl_{n+1}$. In particular, for $\lambda=((k-1)^{n+1})+\nu$
where $|\nu|=n+1$, one has $\K_{\lambda,(k^{n+1})}=\K_{\nu,(1^{n+1})}$,
a `first layer formula' in the sense
of \cite{Ste}. It would be interesting to know whether there exists
multivariate analogs of the other polynomials considered in \cite{Ste},
associated to the exterior algebra of $\sl_{n+1}$ and to the Macdonald
complex.

\smallskip\noindent
{\bf 2.} Putting $x_k=q^k$ in Theorem \ref{RECT}, one recovers a
known formula (\cite{Gupta1}, Prop. 4.2). For the general case,
there is an asymptotic formula of Stanley \cite{Sta84}
\begin{equation}\label{STA}
f_{\alpha\beta}(q) := \lim_{n\rightarrow\infty}
K_{[\alpha,\beta]_{n+1},(\beta_1^{n+1})}(q)
= s_\alpha * s_\beta (q,q^2,q^3,\ldots )
\end{equation}
where $*$ is the internal product of symmetric functions. When
neither $\alpha$ nor $\beta$ is a row,
$\K_{[\alpha,\beta]_{n+1},(\beta_1^{n+1})}$
is not symmetric, and (\ref{STA}) cannot be interpreted as a specialization
of a symmetric identity. Multivariate analogues of the asymptotic
multiplicities
$f_{\alpha\beta}(q)$ have been considered in \cite{Br}.
Other properties of the generalized exponents can be found in \cite{Ki4}.

\end{document}